\documentclass[12pt]{amsart}
\usepackage{amsmath,amsthm,amssymb}
\usepackage{euscript}
\usepackage[all]{xypic}
\usepackage{enumerate}
\usepackage{array}
\DeclareMathAlphabet{\EuRm}{U}{eur}{m}{n}
\SetMathAlphabet{\EuRm}{bold}{U}{eur}{b}{n}
\begin{document}
%
%
\theoremstyle{plain}
\swapnumbers
	\newtheorem{thm}{Theorem}[section]
	\newtheorem{prop}[thm]{Proposition}
	\newtheorem{lemma}[thm]{Lemma}
	\newtheorem{cor}[thm]{Corollary}
	\newtheorem{fact}[thm]{Fact}
\theoremstyle{definition}
	\newtheorem{defn}[thm]{Definition}
	\newtheorem{expl}[thm]{Explanation}
	\newtheorem{claim}[thm]{Claim}
	\newtheorem{notn}[thm]{Notation}
	\newtheorem{construct}[thm]{Construction}
	\newtheorem{subsec}[thm]{}
\theoremstyle{remark}
        \newtheorem{aside}[thm]{Aside}
        \newtheorem{remark}[thm]{Remark}
	\newtheorem{example}[thm]{Example}
	\newtheorem{examples}[thm]{Examples}
        \newtheorem{question}[thm]{Question}
       	\newtheorem{ack}[thm]{Acknowledgements}
\newenvironment{myeq}[1][]
{\stepcounter{thm}\begin{equation}\tag{\thethm}{#1}}
{\end{equation}}
\newcommand{\mydiag}[2][]{\myeq[#1]{\xymatrix{#2}}}
\newcommand{\mydiagram}[2][]
{\stepcounter{thm}\begin{equation}
     \tag{\thethm}{#1}\xymatrix{#2}\end{equation}}
%
\newenvironment{mysubsection}[2][]
{\begin{subsec}\begin{upshape}\begin{bfseries}{#2.}
\end{bfseries}{#1}}
{\end{upshape}\end{subsec}}
\newenvironment{mysubsect}[2][]
{\begin{subsec}\begin{upshape}\begin{bfseries}{#2\vsn.}
\end{bfseries}{#1}}
{\end{upshape}\end{subsec}}
\newcommand{\w}[2][ ]{\ \ensuremath{#2}{#1}\ }
%
%
%
\newcommand{\xra}[1]{\xrightarrow{#1}}
\newcommand{\xla}[1]{\xleftarrow{#1}}
\newcommand{\hra}{\hookrightarrow}
\newcommand{\lora}{\longrightarrow}
\newcommand{\adj}[2]{\hsm\substack{{#1}\\ \rightleftharpoons\\ {#2}}\hsm}
\newcommand{\h}{\hspace{0.22 mm}}
\newcommand{\hsp}{\hspace{10 mm}}
\newcommand{\hs}{\hspace{7 mm}}
\newcommand{\hsm}{\hspace{2 mm}}
\newcommand{\hsmin}{\hspace{-2 mm}}
\newcommand{\vs}{\vspace{7 mm}}
\newcommand{\vsm}{\vspace{3 mm}}
\newcommand{\vsn}{\vspace{1 mm}}
\newcommand{\rest}[1]{\lvert_{#1}}
\newcommand{\lra}[1]{\langle{#1}\rangle}
\newcommand{\EQUIV}{\Leftrightarrow}
\newcommand{\epic}{\to\hspace{-5 mm}\to}
\newcommand{\xepic}[1]{\xrightarrow{#1}\hspace{-5 mm}\to}
\newcommand{\hotimes}{\hat{\otimes}}
\newcommand{\eee}{\hfill$\Box$}
%
%
\newcommand{\ch}{{\EuScript Ch}}
\newcommand{\Ch}[1]{\ch_{#1}}
\newcommand{\cCh}[1]{c\Ch{#1}}
\newcommand{\Cok}{\operatorname{Coker}\,}
\newcommand{\colim}{\operatorname{colim}}
\newcommand{\colimit}[1]
{\colim\hspace{-10mm}\raisebox{-2.0ex}{{\scriptsize ${#1}$}}}
\newcommand{\Der}{\operatorname{Der}}
\newcommand{\diag}{\operatorname{diag}}
\newcommand{\ext}{\operatorname{ex}}
\newcommand{\Ext}{\operatorname{Ext}}
\newcommand{\bExt}{\overline{\Ext}\,}
\newcommand{\cExt}{{\EuScript Ext}\,}
\newcommand{\wExt}{\widehat{Ext}}
\newcommand{\gr}{\operatorname{gr}}
\newcommand{\ho}{\operatorname{ho}}
\newcommand{\holim}{\operatorname{holim}}
\newcommand{\Hom}{\operatorname{Hom}}
\newcommand{\uHom}{\underline{\Hom}}
\newcommand{\Id}{\operatorname{Id}}
\newcommand{\Image}{\operatorname{Im}\h}
\newcommand{\inte}{\operatorname{in}}
\newcommand{\Ker}{\operatorname{Ker}}
\newcommand{\map}{\operatorname{map}}
\newcommand{\maps}{\map_{\ast}}
\newcommand{\Obj}{\operatorname{Obj}}
\newcommand{\Nat}{\operatorname{Nat}}
\newcommand{\op}{^{\operatorname{op}}}
\newcommand{\sk}[1]{\operatorname{sk}_{#1}}
\newcommand{\Spec}{\operatorname{Spec}}
\newcommand{\Tor}{\operatorname{Tor}}
\newcommand{\bTor}{\overline{\Tor}\,}
\newcommand{\cTor}{{\EuScript Tor}\,}
\newcommand{\Tot}{\operatorname{Tot}}
\newcommand{\tru}[1]{\operatorname{tr}_{#1}}
%
%
\newcommand{\bB}{\mathbf{B}}
\newcommand{\bDb}{\mathbf{\Delta}^{\bullet}}
\newcommand{\bp}{\mathbf{\pi}}
\newcommand{\bS}{\mathbf{S}}
\newcommand{\LL}{\mathbb L}
\newcommand{\NN}{\mathbb N}
\newcommand{\QQ}{\mathbb Q}
\newcommand{\RR}{\mathbb R}
\newcommand{\SSa}{\mathbb S}
\newcommand{\ZZ}{\mathbb Z}
%
%
\newcommand{\cA}{{\mathcal A}}
\newcommand{\cB}{{\mathcal B}}
\newcommand{\cC}{{\mathcal C}}
\newcommand{\hC}{\hat{\cC}}
\newcommand{\sC}{s\cC}
\newcommand{\cDD}{{\mathcal D}}
\newcommand{\cE}{{\mathcal E}}
\newcommand{\cF}{{\mathcal F}}
\newcommand{\cG}{{\mathcal G}}
\newcommand{\cI}{{\mathcal I}}
\newcommand{\cM}{{\mathcal M}}
\newcommand{\hM}{\hat{\cM}}
\newcommand{\MP}{\cM_{\Phi}}
\newcommand{\cO}{{\mathcal O}}
\newcommand{\cP}{{\mathcal P}}
\newcommand{\cS}{{\mathcal S}}
\newcommand{\Sa}{\cS_{\ast}}
\newcommand{\cT}{{\mathcal T}}
\newcommand{\Ta}{\cT_{\ast}}
\newcommand{\cU}{{\mathcal U}}
\newcommand{\hU}{\hat{\cU}}
\newcommand{\cV}{{\mathcal V}}
%
%
\newcommand{\hy}[2]{{#1}\text{-}{#2}}
\newcommand{\AToX}{\fA\ToX}
\newcommand{\Ab}{{\operatorname{Ab}}}
\newcommand{\ab}{{\operatorname{ab}}}
\newcommand{\Abgp}{{\EuScript Abgp}}
\newcommand{\Alg}[1]{\hy{#1}{\EuScript Alg}}
\newcommand{\Aa}{A_{\Theta}}
\newcommand{\Ax}{A_{\TX}}
\newcommand{\Aox}{A_{\ToX}}
\newcommand{\hAox}{\hat{A}_{\ToX}}
\newcommand{\Ap}{A_{\Phi}}
\newcommand{\BPA}{\hy{\cB}{\PAlg}}
\newcommand{\CP}{\hy{\cC}{\Pi}}
\newcommand{\CPA}{\hy{\cC}{\PAlg}}
\newcommand{\Gam}[1]{\hy{\Gamma}{#1}}
\newcommand{\GS}{\Gam{\Sa}}
\newcommand{\Gp}{{\EuScript Gp}}
\newcommand{\mt}{m_{\theta}}
\newcommand{\hm}{\hat{m}_{\theta}}
\newcommand{\tii}{\tilde{\text{\i}}}
\newcommand{\tjj}{\tilde{\text{\j}}}
\newcommand{\Phal}{$\Phi$-algebra}
\newcommand{\PhA}{\Alg{\Phi}}
\newcommand{\PAlg}{\Alg{\Pi}}
\newcommand{\PhC}{\hy{\Phi}{\cC}}
\newcommand{\PhE}{\hy{\Phi}{\cE}}
\newcommand{\PhV}{\hy{\Phi}{\cV}}
\newcommand{\PhTA}{\hy{\Phi}{\TA}}
\newcommand{\PhTAX}{\hy{\Phi}{\TA/X}}
\newcommand{\hp}{\hat{\phi}}
\newcommand{\RM}[1]{\hy{#1}{\EuScript Mod}}
\newcommand{\Set}{{\EuScript Set}}
\newcommand{\Seta}{\Set_{\ast}}
\newcommand{\Tab}{\Theta_{\ab}}
\newcommand{\Tg}{\Theta_{\fG}}
\newcommand{\TX}{\Theta_{X}}
\newcommand{\ToX}{\Theta/X}
\newcommand{\TU}{\Theta_{\cU}}
\newcommand{\TUA}{\Alg{\TU}}
\newcommand{\Td}{\Theta^{\delta}}
\newcommand{\TC}{\hy{\Theta}{\cC}}
\newcommand{\TA}{\Alg{\Theta}}
\newcommand{\TtA}{\Alg{\tilde{\Theta}}}
\newcommand{\ATA}{\hy{\fA}{\TA}}
\newcommand{\ATAX}{\hy{\fA}{(\TA/X)}}
\newcommand{\TAb}{\Alg{\Tab}}
\newcommand{\TG}{\Alg{\Tg}}
\newcommand{\TXA}{\Alg{\TX}}
\newcommand{\TAXa}{(\TA/X)_{\ab}}
\newcommand{\ToXA}{\Alg{\ToX}}
\newcommand{\TS}{\hy{\Theta}{\Seta}}
\newcommand{\TV}{\hy{\Theta}{\cV}}
\newcommand{\TdA}{\Alg{\Td}}
\newcommand{\Tal}{$\Theta$-algebra}
\newcommand{\XA}{\Alg{\fX}}
\newcommand{\XC}{\hy{\fX}{\cC}}
%
%
\newcommand{\fa}{\mathfrak{a}}
\newcommand{\fA}{\mathfrak{A}}
\newcommand{\fC}{\mathfrak{C}}
\newcommand{\fF}{{\mathfrak F}}
\newcommand{\fFT}{\fF_{\Theta}}
\newcommand{\fFTp}{\fF'_{\Theta}}
\newcommand{\fg}{\mathfrak{g}}
\newcommand{\fG}{\mathfrak{G}}
\newcommand{\fL}{\mathfrak{L}}
\newcommand{\fM}{\mathfrak{M}}
\newcommand{\fm}{\mathfrak{m}}
\newcommand{\fS}{\mathfrak{S}}
\newcommand{\fW}{\mathfrak{W}}
\newcommand{\fWS}{\mathfrak{WS}}
\newcommand{\fX}{\mathfrak{X}}
%
%
\newcommand{\co}[1]{c({#1})_{\bullet}}
\newcommand{\cu}[1]{c({#1})^{\bullet}}
\newcommand{\Bd}{B_{\bullet}}
\newcommand{\Bdd}{B_{\bullet\bullet}}
\newcommand{\Cu}{C^{\bullet}}
\newcommand{\Ed}{E_{\bullet}}
\newcommand{\Eu}{E^{\bullet}}
\newcommand{\Euu}{E^{\bullet\bullet}}
\newcommand{\Vd}{V_{\bullet}}
\newcommand{\Wd}{W_{\bullet}}
\newcommand{\Wu}{W^{\bullet}}
\newcommand{\Wus}{W^{\ast}}
\newcommand{\Xd}{X_{\bullet}}
\newcommand{\Xu}{X^{\bullet}}
\newcommand{\Yd}{Y_{\bullet}}
\newcommand{\Yu}{Y^{\bullet}}
%
%
\newcommand{\bz}{\mathbf{0}}
\newcommand{\bo}{\mathbf{1}}
\newcommand{\bt}{\mathbf{2}}
\newcommand{\bk}{\mathbf{k}}
\newcommand{\bn}{\mathbf{n}}
%
%
\newcommand{\PiM}{\Pi_{\cM}}
\newcommand{\pin}[2]{\pi_{{#1},{#2}}}
\newcommand{\piM}[1]{\pin{M}{#1}}
\newcommand{\picM}[1]{\pin{\cM}{#1}}
\newcommand{\q}[1]{^{({#1})}}
%
%
\newcommand{\FT}{F_{\Theta}}
\newcommand{\FTab}{F_{\Tab}}
\newcommand{\FTx}{F_{\TX}}
\newcommand{\FTox}{F_{\ToX}}
\newcommand{\UT}{U_{\Theta}}
\newcommand{\UTab}{U_{\Tab}}
\newcommand{\UTx}{U_{\TX}}
%
%
\title{Generalized Andr\'{e}-Quillen Cohomology}
\author{David Blanc}
\date{August 10, 2007; revised: February 10, 2008 }
\address{Dept.\ of Mathematics, Faculty of Sciences, Univ.\ of Haifa, 
31905 Haifa, Israel}
\email{blanc@math.haifa.ac.il}

\begin{abstract}
We explain how the approach of Andr\'{e} and Quillen to defining
cohomology and homology as suitable derived functors extends to generalized
(co)homology theories, and how this identification may be used to
study the relationship between them.
\end{abstract}

\maketitle 

\section*{Introduction}

After the cohomology of topological spaces was discovered in the 1930's, 
the concept was expanded to groups, and later to associative, commutative, and 
Lie algebras, in the 1940's and early 1950's. In the following decade the first 
generalized cohomology theories for spaces appeared (see \cite{MacLO,MasHCT}). 
All these examples started out in the form of explicit constructions, and only 
later were their theoretical underpinnings provided: in particular, cohomology
for general algebraic categories was described by Beck and others
in terms of triples (see \cite{BecT}, and compare \cite{DuskS}), and
then by Andr\'{e} and Quillen in terms of (non-abelian) derived
functors (see \cite{AndrM,QuiH}). In the latter version, cohomology
groups are the derived functors of \ $\Hom$ \ into a fixed abelian
group object (and homology groups are the derived functors of abelianization).

However, for topological spaces the only abelian group objects are
(products of) Eilen\-berg-Mac~Lane spaces, which represent ordinary
cohomology. Thus we need a different framework to describe generalized 
(co)homology: this is provided by stable homotopy theory (cf.\ \cite{BrowC,GWhG}).

Our goal here is to provide a uniform definition for homology and
cohomology encompassing the theories mentioned above, as well as some
new ones. As a side benefit, we clarify exactly what assumptions on an
(algebraic) category $\cC$ are needed in order for the approach of
Andr\'{e} and Quillen to work. (This is the reason for the somewhat 
technical Section \ref{caft}.)

The approach given here applies, inter alia, to:

\begin{enumerate}\renewcommand{\labelenumi}{(\alph{enumi})}
\item Homology and cohomology of groups and various types of algebras;
\item Versions of the above with local coefficients (\S \ref{scoht}-\ref{shomt});
\item Unstable generalized (co)homology of spaces (\S \ref{sgch}-\ref{sgh}); 
\item Generalized (co)homology of spectra and spaces (\S \ref{echsm});
\item Cohomology of operads, and of algebras over an operad (\S \ref{scoa});
\item Cohomology of diagrams of spaces or algebras (\S \ref{sdiagt}).
\end{enumerate}

The last two have applications to deformation theory (see \cite{MarkCT,MShniCD} 
and \cite{GSchaD,GGSchaD}, respectively). 

The cohomology of sheaves has a dual definition to the one presented here
here (see \S \ref{scohsh}). Of course, there are other concepts of
cohomology which do not fit into our framework; most notably, a number
of versions of the cohomology of categories (see \S \ref{rjpir}). 

\begin{mysubsection}{Representing cohomology}
\label{src}
In order to define a cohomology theory in a category $\cC$, we need 
a representing object \ $G\in\cC$, \ as well as a suitable model
category structure on the category  \ $\sC=\cC^{\Delta\op}$ \ of
simplicial objects over $\cC$ (see \S \ref{srmc}). However, in this generality \
$\Hom\h_{\cC}(-,G)$ \ will take values in sets, and applying this
functor to a simplicial resolution \ $\Vd\to X$ \ in \ $\sC$ \ just
yields a cosimplicial set, for which we have no appropriate model
category. It turns out that in order to get an interesting cohomology
theory, two ingredients are generally needed: 
\end{mysubsection}

\begin{enumerate}\renewcommand{\labelenumi}{$\bullet$ \ }
\item The category $\cC$ must be enriched over a symmetric monoidal
  category $\cV$;  
\item The representing object $G$ must have additional ``algebraic''
  structure.
\end{enumerate}

We shall use the concept of a \emph{sketch} \ -- \ a straightforward
generalization of Lawvere's concept of a \emph{theory} \ -- \ to
describe this additional structure (see \S \ref{dco}). In this language,
we say that $G$ is a \Phal\ in $\cC$, for a suitable FP-sketch
$\Phi$. We also use sketches to describe the kind of algebraic
categories to which our approach applies: this will allow us to treat
operads and their algebras, for example, uniformly with the usual
universal algebras.   

\begin{enumerate}\renewcommand{\labelenumi}{$\bullet$ \ }
\item Note that the functor \ $\Hom_{\cC}(-,G)$ \ now takes values in
  the category $\cDD$ of (cosimplicial) \Phal s in $\cV$. Our final
  requirement is that the above two ingredients must combine to make 
  $\cDD$ into a (semi-) triangulated model category (see \S \ref{dqtc}). 
\end{enumerate}

The question we consider here is in some sense dual to that of Brown
Representability in triangulated categories (cf.\ 
\cite{CKNeemF,FrankB,KrausB,NeemTC}): 
rather than asking which cohomology functors are representable, we
seek conditions for a representable functor to be a cohomology theory.

\begin{examples}
\label{ecc}
In the category of groups (where \ $\cV=\Set$), \ with an abelian 
group $G$ as the coefficients, the model category we consider is that 
of simplicial groups. The total left derived functor of \ $\Hom(-,G)$ \ 
then takes values in the semi-triangulated category of cosimplicial
abelian groups (equivalently, cochain complexes). 

On the other hand, for pointed simplicial sets or topological spaces 
(where \ $\cV=\Sa$), \ we may take \ $\Phi=\Gamma$, \ and  \ 
$\Hom(-,G)$ \ takes values in $\Gamma$-spaces \ -- \ again, a
semi-triangulated category.

Note that the category of spectra is triangulated (and enriched over
itself), so we can take \emph{any} spectrum $G$ as coefficients. 
\end{examples}

Our original motivation for creating a joint setting for algebraic and
generalized topological (co)homology theories was to try to gain a
better understanding of the relationship between homology and 
cohomology. This is provided by a universal coefficients spectral
sequence (see Theorem \ref{tone} below). We obtain a similar result
for homology (Proposition \ref{ptwo}), as well as ``reverse Adams
spectral sequences'' (Theorems \ref{ttwo} and \ref{tthree}) relating
homotopy to (co)homology.  

\begin{mysubsection}{Notation and conventions}
\label{snac}
The category of topological spaces is denoted by $\cT$, and that of 
pointed connected topological spaces by  \ $\Ta$. \  
The category of groups is denoted by \ $\Gp$, \ that of abelian groups
by \ $\Abgp$, \ and that of pointed sets by \ $\Seta$. \ 
For any category $\cC$, \ $\gr_{S}\cC$ \ denotes the category of
$S$-graded objects over $\cC$  (i.e., diagrams indexed by the discrete
category $S$), \ $\sC$ \ that of simplicial objects over $\cC$, and \
$c\cC$ \ that of cosimplicial objects over $\cC$. The category of simplical
sets will be denoted by $\cS$, that of reduced simplicial sets by \
$\Sa$, \  and that of simplicial groups by $\cG$. 
For any \ $Z\in\cC$, \ we write \ $\co{Z}$ \ for the constant
simplicial object determined by $Z$, and \ $\cu{Z}$ \ for the constant
cosimplicial object. \ If $\cA$ is any abelian category, we denote
the category of chain complexes over $\cA$ by \ $\ch(\cA)$; \ however,
we write \ $\Ch{R}$ \ for \ $\ch(\RM{R})$, \ and similarly \ $\cCh{R}$ \ 
for cochain complexes of $R$-modules.
\end{mysubsection}

\begin{mysubsection}{Organization}\label{sorg}
Section \ref{cat} provides background material on sketches, theories,
and algebras over them. In Section \ref{cgch} we give our abstract
definition of homology and cohomology, in the context of suitable
model categories. Abelian group objects in sketchable categories are
described in Section \ref{caft}, and these are used in Section
\ref{cch} to define the (co)homology of \Tal s. Section \ref{cgencoh}
explains how generalized cohomologies fit into our framework, using
$\Gamma$-spaces. Finally, the theory is applied in Section \ref{css}
to construct universal coefficient and reverse Adams spectral
sequences in this general framework. 
\end{mysubsection}

\begin{ack}
This paper is an outgrowth of joint work with George Peschke, in
\cite{BPescF}, and I would like to thank him for many useful
discussions and insights. I also thank the referee for his or her
helpful comments, and the Institut Mittag-Leffler (Djursholm, Sweden) 
for its hospitality during the period when this paper was completed.
\end{ack}

%
%
\section{Algebras and theories}
\label{cat}

As Lawvere observed (cf.\ \cite{LawvF}), `varieties of universal 
algebras' in the sense of Mac~Lane (cf.\ \cite[V,6]{MacLC}) can be 
corepresented by functors out of a fixed category $\Theta$. 
This idea was later generalized by Ehresmann to sketches (see
\cite{BEhreC}), which turn out to be the most convenient language to
describe both the algebraic categories we work in, and the
representing objects for cohomology.

\begin{defn}\label{dco}
A \emph{sketch} \ $\lra{\Theta,\cP,\cI}$ \ is a small category
$\Theta$ with distinguished sets $\cP$ of (limit) cones and $\cI$ of 
(colimit) cocones. In particular, a \emph{finite product (FP-)sketch}
is a sketch in which $\cP$ consists only of finite products (and \ 
$\cI=\emptyset$). \ A \emph{theory} is an FP-sketch $\Theta$
containing a zero object, for which $\cP$ consists of \emph{all}
finite products.

We think of a map \ $f:\vartheta_{1}\times\dotsc\times\vartheta_{n}\to\theta$ \ 
in $\Theta$ as corepresenting a (possibly graded) $n$-ary operation.
A theory $\Theta$ is \emph{sorted}  by a set \ $S\subseteq\Obj\Theta$ \ 
if every object in $\Theta$ is uniquely isomorphic to a finite
product of objects from $S$ (see \cite[\S 5.6]{BorcH2}).   
Lawvere originally considered only theories sorted by \ $\{\bo\}$, \
so that \ $\Obj(\Theta)=\NN$, \ with \  $\bn\cong\prod_{i=1}^{n}~\bo$ \
for \ $n\geq 0$. 

If $\Theta$ is an FP-sketch and $\cC$ is any pointed category, \
a \emph{\Tal\ in $\cC$} is a pointed functor \ 
$X:\Theta\to\cC$ \ which preserves all products in $\cP$. More generally, if
$\Theta$ is any sketch, a \Tal\ \ $X:\Theta\to\cC$ \ is required to
preserve all distinguished limits (in $\cP$) and colimits (in $\cI$). 
The category of \Tal s in $\cC$ is denoted by \ $\TC$; \ a \Tal\ in \ $\Seta$ \ 
will be called simply a \emph{\Tal}, and we write \ $\TA$ \ for \ 
$\hy{\Theta}{\Seta}$. \ We call a category $\cDD$ \emph{sketchable} if
it is equivalent to \ $\TA$, \ and say that $\Theta$ \emph{sketches}
$\cDD$. Such categories are accessible, in the sense of model theory, as well
as being locally presentable  (see \cite[Cor.~2.61 \& 1.52]{ARosiL}).   
A \emph{map of theories} (or of sketches) \ $\psi:\Theta\to\Theta'$ \ 
is a functor which preserves all products (respectively, all
distinguished limits and colimits). Such a map $\psi$ induces a functor \
$\psi^{\ast}:\Alg{\Theta'}\to\TA$.

More generally, if $\Theta$ is a theory (or FP-sketch),  a
\emph{\Tal} in any symmetric monoidal category \ $\lra{\cV,\otimes,I}$ \ 
(cf.\ \cite[\S 6.1]{BorcH2}) is a functor \ $X:\Theta\to\cV$ \ taking
the (distinguished) products in $\Theta$ to $\otimes$-products in
$\cV$, with \ $X(\ast)=I$.
\end{defn}

\begin{remark}\label{rhom}
Since we can think of a \Tal\ $X$ in $\cC$ as a certain kind of diagram
in $\cC$ (with specified products), we see that \ $\Hom\h_{\cC}(-,X)$ \
takes values in \ $\TA$. \ More generally, if $\cC$ is enriched over a
symmetric monoidal category  \ $\lra{\cV,\otimes,I}$ \ via \
$\map_{\cC}$ \ (cf.\ \cite[\S 6.2]{BorcH2}), and  \ $\map_{\cC}(A,-)$ \ 
takes products to $\otimes$, then \ $\map_{\cC}(-,X)$ \ take values in \ $\TV$.  
\end{remark}

\begin{examples}\label{egtheory}

\noindent (a)\ The category of groups is sketched by a theory $\fG$,
with \ $\mu:\bt\to\bo$ \ representing the group operation, \
$\rho:\bo\to\bo$ \ the inverse, and \ $e:\bz\to\bo$ \ the identity
(satisfying the obvious relations).  Similarly, the category of
abelian groups is sketched by a theory \ $\fA$ \ (with the same maps,
satisfying a further relation) and the inclusion \ $i:\fG\subset\fA$ \
induces the inclusion of categories \ $\Abgp\subseteq\Gp$\vsm.

\noindent(b)\  An operad \ $\Gamma=(\Gamma(n))_{n=0}^{\infty}$ \ is an 
$\cO$-algebra in a symmetric monoidal category \
$\lra{\cV,\otimes,I}$, \ where \ $\cO$ is a ``universal'' theory for
operads. Similarly, an algebra over the operad $\Gamma$ (see 
\cite[\S 14]{MayG})  is just a \ $\Theta_{\Gamma}$-algebra in \
$\lra{\cV,\otimes,I}$, \ where the theory \ $\Theta_{\Gamma}$ \ is
obtained from $\Gamma$ in the obvious way (replacing $\otimes$ with
$\times$). The same applies more generally to PROP's, colored operads,
and other variants (see \cite{MSStasO} for a survey on operads,
especially in the algebraic context)\vsm. 

\noindent (c)\ Given a topological space $X$, let $\cU$ denote the
directed set of non-empty open sets in $X$, with inclusions \ -- \
so that \ $\cU\op$ \ sketches presheaves of sets.
By adding arbitrary formal coproducts \ 
$\coprod_{\alpha\in A}\,U_{\alpha}$ \ for any collection \
$\{U_{\alpha}\}_{\alpha\in A}$ \ in $\cU$, we obtain a category $\hU$,
in which the diagram; 
%
\begin{myeq}\label{etwelve}
\xymatrix@R=25pt{
\coprod_{(\alpha,\beta)\in A\times A}~
U_{\alpha}\cap U_{\beta} \ar@<1mm>[r]^<<<<{i} \ar@<-1mm>[r]_<<<<{j} &  
\coprod_{\alpha\in A}~U_{\alpha} \ar[r]^{\kappa} &
\bigcup_{\alpha\in A}~U_{\alpha} 
}
\end{myeq}
\noindent is a coequalizer (if the first term is empty, $\kappa$
is an isomorphism). 

If we now let \ $\TU:=\hU\op$ \ (sorted by $\cU$), with $\cP$
consisting of the opposites of the formal coproducts and of all the
coequalizers \ \eqref{etwelve} \ (and \ $\cI=\emptyset$), \  
we obtain a sketch whose algebras \ $\cF:\Theta_{\cU}\to\Set$ \ 
are sheaves of sets on $X$. Furthermore, for any \ $V\in\cU$, \ if:
$$
C_{V}(U)~:=~\begin{cases}\{\ast\}& \ \text{if} \ U\subseteq V\\
~\emptyset& \ \text{if} \ U\not\subseteq V,\end{cases}
$$
\noindent there is a natural isomorphism \ $\Hom\h_{\TUA}(C_{V},\cF)=\cF(V)$.
\end{examples}

\begin{defn}\label{dgpthy}
Given a theory $\fX$, an $\fX$-\emph{theory} (or sketch) $\Theta$
is one equipped with a map of theories (or sketches) \ 
$\psi:\coprod_{S}\,\fX\to\Theta$ \ which is bijective on objects,
where the coproduct is taken in the category of theories (or sketches) 
over some index set $S$.  If $\fX$ is sorted by $\{\bo\}$, an
\emph{$\fX$-structure at} an object $c$ in a category $\cC$ is an
$\fX$-algebra \ $\rho:\fX\to\cC$ \ with \ $\rho(\bo)=c$. \ A theory
$\Theta$ sorted by $S$ is an $\fX$-theory if and only if it is
equipped with an $\fX$-structure at each \ $s\in S$. 

If all other maps of $\Theta$ commute with those coming from $\psi$,
we call $\Theta$ a \emph{strong $\fX$-theory} (or sketch).
\end{defn}

\begin{example}\label{eggth}
If $\Theta$ is a $\fG$-theory, then the map of theories \
$\psi:\coprod_{S}\,\fG\to\Theta$ \ induces an ``underlying $S$-graded
group'' functor \ $\psi^{\ast}$, \ which we denote by \ 
$V:\TA\to\Gp^{S}=\Alg{\coprod_{S}\,\fG}$. \ $\Theta$ is a strong
$\fG$-theory if all the operations in $\Theta$ are homomorphisms of
the underlying graded group. 
\end{example}

\begin{mysubsection}{Free $\Theta$-algebras}
\label{sfta}
For any theory $\Theta$, let \ $\Td$ \ denote the discrete theory
with the same objects (and products) as $\Theta$. If $\Theta$
is sorted by $S$, \ $\Td$ \ sketches the category of $S$-graded
sets, and the inclusion \ $I:\Td\hra\Theta$ \ induces the
forgetful functor \ $U=\UT:\TA\to\TdA$. \ As usual, there is a
\emph{free} functor \ $F=\FT:\TdA\to\TA$ \ left adjoint to \ $\UT$. \
We denote by \ $\fFT$ \ the full subcategory of \ $\TA$ \ whose
objects are free (that is, in the image of \ $\FT$).

Since all limit-sketchable categories are locally presentable,
they are complete (see, e.g., \cite[Theorem~1.46]{ARosiL}) and
cocomplete. Thus for any theory $\Theta$, the category \ $\TA$ \ of
$\Theta$-algebras has all limits and colimits.
\end{mysubsection}

\begin{mysubsection}{Sketching \Phal s in \ $\TA$}
\label{sspat}
If $\Theta$ is a theory (sorted by $S$) and $\Phi$ is another
theory (singly sorted, for simplicity), the category \ $\PhTA$ \ of
\Phal s in \ $\TA$ \ is sketched by a theory \ $\Phi(\Theta)$ \ 
(sorted by $S$), defined as follows:

\begin{enumerate}\renewcommand{\labelenumi}{(\alph{enumi})}
\item We first add an $S$-graded copy of $\Phi$ to $\Theta$, \ 
setting \ $\Theta_{\Phi}:=\Theta\cup_{S}\coprod_{S}\Phi$, \ 
so that we now have each operation of $\Phi$ acting on each \
$\theta\in S$. \  The inclusion \ $i:\Theta\hra\Theta_{\Phi}$ \
induces a forgetful functor \ $i^{\ast}:\Alg{\Theta_{\Phi}}\to\TA$.
\item Next, we force all operations of $\Theta$ to commute with the 
new operations \ - \ that is, for each \ $f:\theta_{1}\to\theta_{2}$ \ 
in $\Theta$ and \ $g:\bn\to\bk$ \ in $\Phi$, we require that
$$
\xymatrix@R=25pt{
\theta_{1}^{n} \ar[r]^{g} \ar[d]^{f^{n}} & \theta_{1}^{k} \ar[d]^{f^{k}}\\ 
\theta_{2}^{n} \ar[r]^{g} & \theta_{2}^{k} 
}
$$
commute, so we obtain a quotient of theories \ 
$q:\Theta_{\Phi}\epic\Phi(\Theta)$.
\end{enumerate}

By construction \ $\Alg{\Phi(\Theta)}\cong\PhTA$. \ Note that \ 
$q^{\ast}$ \ and \ $i^{\ast}$ \ commute with the underlying $S$-graded
set functors \ $\UT$, \ $U_{\Theta_{\Phi}}$, and \ $U_{\Phi\Theta}$, \
which create all limits in their respective categories, so \ 
$q^{\ast}$ \ and \ $i^{\ast}$ \ commute with all (small) limits. Thus
by \cite[Theorem~5.5.7]{BorcH2} each has a left adjoint. The adjoint of the 
composite \ $i^{\ast}\circ q^{\ast}:\PhTA\to\TA$ \ will be called the
\emph{$\Phi$-localization} of \ $\TA$, \ and denoted by \ 
$L_{\Phi}:\TA\to\PhTA$.
\end{mysubsection}

\begin{remark}\label{renrich}
Note that given $G$ in \ $\PhTA$, \ by Remark \ref{rhom} \
$\Hom_{\TA}(-,G)$ \ has a natural structure of a \Phal. Furthermore,
if \ $i^{\ast}\circ q^{\ast}$ \ is a faithful embedding of categories
(which will happen if $\Theta$ is a $\Phi$-theory, for example), then \ 
$L_{\Phi}$ \ is idempotent and any \Phal\ in \ $\TA$ \ is in the image
of \ $L_{\Phi}$, \ up to natural isomorphism. Thus \ $\Hom_{\PhTA}(-,-)$ \ 
has a natural structure of a \Phal, in this case. By mimicking the
construction of \  $A\times B\to A\otimes B$ \ for abelian groups, one
can then make \ $\PhTA$ \ into a closed symmetric monoidal category
(see \cite[\S 6.1.3]{BorcH2}).
\end{remark}

%
%
\section{Generalized homology and cohomology}
\label{cgch}
 
We are now able to give a definition of homology and cohomology for
model categories, somewhat more general than Quillen's original
approach (cf.\ \cite[II, \S 5]{QuiH}):

\begin{mysubsect}{Triangulated categories}
\label{striang}

The target of a cohomology functor should be a model category whose
homotopy category is triangulated. There are a number of variants of
this concept, originally due to Grothendieck.   
For our purposes, a \emph{triangulated category} is an additive
category $\cC$ equipped with an automorphism \ $T:\cC\to\cC$ \ (called
the \emph{translation} functor), and a collection $\cDD$ of
\emph{distinguished triangles} of the form \ 
$\lra{X\xra{f} Y\xra{g} Z\xra{h} TX}$, \ satisfying the four axioms of 
\cite[\S 1]{HartsR} (which codify the properties of cofibration
sequences in pointed model categories \ -- \ see \cite[I, \S 3]{QuiH}). 
\end{mysubsect}

\begin{defn}\label{dqtc}
A \emph{semi-triangulated} category is an additive category $\hC$ 
equipped with a collection $\cDD$ of \emph{distinguished triangles} of
satisfying the above four axioms, as well as a translation functor  \
$T:\hC\to\hC$ \ which is an isomorphism onto its image. In all cases
of interest, $T$ can be formally inverted to yield a full triangulated
category \ $\cC=\hC[T^{-1}]$ \ with $\hC$ as a full subcategory;
however, this property is not needed in what follows.  

A set $\cP$ of cogroup objects in $\hC$ will be called a 
\emph{set of generators} if the collection of functors \
$\{\Hom\h_{\hC}(T^{i}P,-)\}_{P\in\cP,i\geq 0}$ \ detects all
isomorphisms in $\hC$.
\end{defn}

\begin{example}\label{esgps}
Typically, (semi-)triangulated categories appear as the homotopy
category of a suitable (semi-)stable model category, as defined
axiomatically in \cite[Ch.\ 7]{HovM} (see also \cite{HPStrA}).  
Thus, the motivating example of a triangulated category is the
homotopy category of (unbounded) chain complexes over an abelian
category $\cA$. Another example is provided by Boardman's stable
homotopy category \ $\ho\Spec$ \ (cf.\ \cite{VogtB}), where there are
a number of different underlying stable model categories (see
\cite{HSSmiS}, \cite{SchweSH}, or \cite{EKMMayR}).

The subcategory $\hC$ of non-negatively graded chain complexes is
semi-triangulated; if $\cA$ has a projective generator $P$, then \ 
$K(P,0)$ \ (the chain complex with $P$ concentrated in degree $0$) is
a generator for $\hC$. 

Similarly, the homotopy category of connective spectra, \
$\ho\Spec_{(0)}$, \ is semi-\-triangu\-la\-ted (with generator \ $S^{0}$).
\end{example}
 
\begin{mysubsection}{Cohomology}
\label{scohf}
In order to define cohomology functors on a model category $\cE$, we assume
that $\cE$ is equipped with: 

\begin{enumerate}\renewcommand{\labelenumi}{(\alph{enumi})~}
\item An FP-sketch  $\Phi$ and a category $\cV$ such that $\cV$ and \ 
 $\PhV$ \ are symmetric monoidal, \ $\cE$ \ is enriched over $\cV$ via \ 
 $\map_{\cE}(-,-):\cE\op\times\cE\to\cV$, \ and \ $\PhE$ \ is enriched 
 over \ $\PhV$ \ via \ $\uHom(-,-):(\PhE)\op\times\PhE\to\PhV$. 
\item An FP-sketch  $\Phi$ and a model category structure on \ $\PhV$ \
  for which \ $\ho\PhV$ \ is semi-triangulated. 
\end{enumerate}

Then for any \ $G\in\PhE$, \ we define the \emph{cohomology of} \
$X\in\cE$ \ \emph{with coefficients in} $G$ to be the total left derived
functor \ $\LL\map_{\cE}(-,G)$ \ of \ $\map_{\cE}(-,G)$, \ applied to
$X$. Recall that \emph{total left derived} functor of a ``left exact''
functor \ $F:\cC\to\cDD$ \ between model categories is defined by
applying $F$ to a cofibrant replacement of $X$ (see \cite[I, \S 4]{QuiH} or 
\cite[8.4]{PHirM}).

If \ $\ho\PhE$ \ has a set of generators $\cP$, then the $\cP$-graded
group \ $H^{n}(X;G):=[T^{n}P,(\LL\map_{\cE}(-,G)) X]_{P\in\cP}$ \ is
called the \emph{$n$-th cohomology group of $X$} \ with coefficients in $G$. 
\end{mysubsection}

\begin{mysubsection}{Homology}
\label{shomol}
To define homology, we need also a homotopy functor \ 
$\Ap:\cE\to\PhE$ \ equipped with a natural isomorphism \ 
$\map_{\cE}(E,X)\xra{\cong}\uHom(\Ap E,X)$ \ in \ $\PhV$ \ 
(cf.\ \S \ref{rhom}) for \ $E\in\cE$ \ and \ $X\in\PhE$. \ 
We then define the \emph{homology} of \ $X\in\cE$ \ to be the total
left derived functor of \ $\Ap$, \ applied to $X$ (\S
\ref{scohf}). Again the \emph{$n$-th homology group of $X$} is: 
$$
H_{n}X~:=~[T^{n}\Ap P,\,(\LL\Ap) X]_{P\in\cP}~.
$$

If \ $\PhE$ \ is a symmetric monoidal model category (see 
\cite[\S 4.2.6]{HovM}), with \ $\uHom(-,Y)$ \ right adjoint
(over \ $\PhV$) \ to \ $-\otimes Y$, \ then for any \ $G\in\PhE$, \  
\emph{homology with coefficients in $G$} is the total left derived
functor of \ $\Ap(-)\otimes G$ \ (assuming $\Ap E$ \ is always
cofibrant). The homology groups \ $H_{n}(X;G)$ \ are defined as above.
Compare \cite[I]{BBeckH}.
\end{mysubsection}

\begin{example}\label{ecohom}
If \ $\cE=\cV=\Sa$ \ (or \ $\Ta$) and \ $\Phi=\fA$, \ then \ 
$\PhC\cong\PhV\cong s\Abgp$ \ and $G$ is a (generalized) Eilenberg-Mac~Lane 
space, so we have ordinary cohomology.
The functor \ $\Ap:\cE\to\PhC$ \ is the usual `abelianization' \
$X\mapsto\ZZ X$, \ which yields ordinary (singular) homology.
\end{example}

\begin{mysubsect}{Resolution model categories}
\label{srmc}

To provide a uniform treatment of the various kinds of (co)homology it
will be convenient to use a framework originally
conceived by Dwyer, Kan and Stover in \cite{DKStoE} under the name of
``$E^{2}$ model categories'', and later generalized by Bousfield (see
\cite{BouC,JardB}. 

Recall that the concept of a \emph{model category} was introduced by
Quillen in \cite{QuiH} to allow application of the methods and
constructions of homotopy theory (of topological spaces) in more
general contexts. This is a category $\cC$, equipped with three
distinguished classes of morphisms \ -- \ weak equivalences, cofibrations,
and fibrations \ -- \  satisfying certain axioms (analogous to those
which hold for the corresponding classes in $\cT$). See \cite{PHirM} or
\cite{HovM} for further details. 

Let $\cC$ be a pointed cofibrantly generated right proper model
category (cf.\ \cite[7.1, 11.1]{PHirM}), equipped with 
a set $\cM$ of cofibrant homotopy cogroup objects
in $\cC$, called \emph{models} (playing the role of the spheres in \
$\Ta$). \ Let \ $\PiM$ \ denote the smallest full subcategory of $\cC$
containing $\cM$ and closed under coproducts, and suspensions (cf.\ 
\cite[I, \S 3]{QuiH}). \ For any \ $X\in\cC$, \ $M\in\cM$, \
and \ $k\geq 0$, \ set \  $\piM{k}X:=[\Sigma^{k}M,X']$, \ where \
$X\to X'$ \ is a fibrant replacement. \ We write \ $\picM{k}X$ \ for
the $\cM$-graded group \ $(\piM{k}X)_{M\in\cM}$.
\end{mysubsect}

\begin{defn}\label{dffm}
A map \ $f:V\to Y$ \ in \ $\sC$ \ is \emph{homotopically $\cM$-free} 
if for each \ $n\geq 0$, \ there is:
\begin{enumerate}
\renewcommand{\labelenumi}{\alph{enumi})\ }
\item a cofibrant object \ $W_{n}\in\PiM$, \ and
\item a map \ $\varphi_{n}:W_{n}\to Y_{n}$ \ in $\cC$ inducing a trivial 
cofibration \ $(V_{n}\amalg_{L_{n}V}L_{n}Y)\amalg W_{n}\to Y_{n}$, \ 
where the $n$-th \emph{latching object} for $Y$ is \ 
$L_{n}Y:=\coprod_{0\leq i\leq n-1} Y_{n-1}/\sim$, \ with \ 
$s_{j_{1}}s_{j_{2}}\dotsc s_{j_{k}}x\in(Y_{n-1})_{i}$ \ is
equivalent to \ $s_{i_{1}}s_{i_{2}}\dotsc s_{i_{k}}x\in(Y_{n-1})_{j}$ \ 
whenever \ 
$s_{i}s_{j_{1}}s_{j_{2}}\dotsc s_{j_{k}}=s_{j}s_{i_{1}}s_{i_{2}}\dotsc
s_{i_{k}}$.
\end{enumerate}

The \emph{resolution model category structure} on \ $\sC$ \ determined
by $\cM$ is now defined by declaring a map \ $f:X\to Y$ \ to be:

\begin{enumerate}
\renewcommand{\labelenumi}{(\roman{enumi})}
\item a \emph{weak equivalence} if \ $\picM{k} f$ \ is a weak equivalence 
of $\cM$-graded simplicial groups for each \ $k\geq 0$;
\item a \emph{cofibration} if it is a retract of a homotopically $\cM$-free map;
\item a \emph{fibration} if it is a Reedy fibration (cf.\ \cite[15.3]{PHirM})
and \ $\piM{k}f$ \ is a fibration of simplicial groups for each \
$M\in\cM$ \ and \ $k\geq 0$.
\end{enumerate}
\end{defn}

\begin{remark}\label{rsrmc}
The resolution model category \ $\sC$ \ is simplicial (cf.\ 
\cite[II,~\S 1]{QuiH}, and is itself endowed with a set of models, 
of the form \ $\hM:=\{S^{n}\otimes M~|\ M\in\cM,\ n\in\NN\}$, \ 
where \ $S^{n}\in\cS$ \ is the simplicial sphere.
\end{remark}

\begin{examples}\label{ermc}
Typical resolution model categories include the following:

\begin{enumerate}
\renewcommand{\labelenumi}{(\roman{enumi})~}
\item When \ $\cC=\Gp$, \ let \ $\cM:=\{\ZZ\}$, \ so \ $\PiM$ \ is
the subcategory of all free groups. The resulting resolution model category
structure on the category \ $\cG=s\Gp$ \ of simplicial groups is the
usual one (see \cite[II,~\S 3]{QuiH}). 
\item More generally, if $\Theta$ is a $\fG$-theory (\S \ref{dgpthy}), 
let \ $\cM:=\fFTp$ \ denote the collection of all monogenic free \Tal s \
$\FT(s)$ \ in \ $\fFT$, \ with $s$ a singleton in \ $\TdA$ \ (i.e., a
graded set, indexed by the discrete sketch \ $\Td$, \ consisting of a
single element in some degree). \ In this case \ $\PiM\cong\fFT$, \
and the model category on \ $s\TA$ \ is that of \cite[II,~\S 4]{QuiH}).   
\item For \ $\cC=\Ta$, \ let \ $\cM:=\{S^{1}\}$, \ so that \ $\PiM$ \
is the homotopy category of wedges of spheres. In this case the model
category of simplicial spaces is the original \ $E^{2}$-model category
of Dwyer, Kan and Stover (cf.\ \cite{DKStoE}).
\end{enumerate}
\end{examples}

\begin{remark}\label{rcomma}
The above discussion is also valid if we work in the comma category \ 
$\TA/X$ \ (cf.\ \cite[II,6]{MacLC}), \ for a $\fG$-theory $\Theta$ and
some fixed \Tal\ $X$.
In fact, any \ $p:\FT\to X$ \ in \ $\fFT/X$ \ is determined by its
adjoint \ $\tilde{p}:T\to\UT X$ \ -- \ in other words, by the \ $\UT
X$-graded set \ $\{p^{-1}(x)]\}_{x\in \UT X}$. \ 
Therefore, \ $\TA/X$ \ can be sketched by a theory \ $\ToX$, \ sorted
by \ $\UT X=\{\phi_{x}~|\ x\in\UT X\}$. \ Note that \ $\ToX$ \ is a 
\emph{$\fG$-sketch over $X$} in the sense that it has $\fG$-structures
of the form:
$$
m_{(x_{1},x_{2})}:\phi_{x_{1}}\times\phi_{x_{2}}\to\phi_{m_{\theta}(x_{1},x_{2})}
$$
for every \ $\theta\in\Theta$ \ and \ $x_{1},x_{2}\in \UT X\theta$ \ 
(and similarly for other morphisms in $\Theta$).

Equivalently, we can equate the discrete theory \ $\ToX^{\delta}$ \
with \ $\TdA/\UT X$, \ and use the adjointness of \ $(\FT,\UT)$ \ to 
define an adjoint pair:
$$   
\ToXA~=~\TA/X\adj{\FT}{\UT}\TdA/\UT X~=~\ToX^{\delta}~.
$$
We can then take the monogenic free \Tal s \ $\fFTp/X$ \ (cf.\ 
\S \ref{ermc}(ii)) \ as our models, and obtain a resolution model
category structure on \ $s(\TA/X)$. \ In particular, any free
resolution \ $\Vd\to X$ \ in \ $s\TA$ \ is also a resolution
(cofibrant replacement) in \ $s(\TA/X)$.
\end{remark}

\begin{mysubsect}{A simplicial version of (co)homology}
\label{scver}

In order to make the abstract description of (co)homology given in 
\S \ref{scohf}-\ref{shomol} more concrete, it is convenient to
formalize the ingredients needed in the following: 
\end{mysubsect}

\begin{defn}\label{dcset}
A \emph{cohomological setting} \ $\lra{\cC,\cM,\cV,\Phi,\Ap}$ \ consists of:
\begin{enumerate}
\item A model category $\cC$, enriched via \ $\map_{\cC}(-,-)$ \ 
over a symmetric monoidal category $\cV$.
\item A set of models $\cM$ for $\cC$.
\item An FP-sketch $\Phi$, such that\vsn:

\begin{enumerate}
\renewcommand{\labelenumii}{(\roman{enumii})~}
\item $\lra{\PhC,\otimes,I,\uHom}$ \ is a closed symmetric monoidal 
     category (with \ $\uHom(G,-)$ \  right adjoint to \ $-\otimes G$).
\item $c\PhV$ \ has a model category structure for which \ $\ho c\PhV$ \ 
   semi-triangulated\vsn. 
\end{enumerate}
\item A homotopy functor \ $\Ap:\PiM\to\PhC$, \ equipped with a
  natural isomorphism: 
%
\begin{myeq}\label{eone}
\nu:\map_{\cC}(F,G)\cong\uHom(\Ap F,X)
\end{myeq}
\noindent for \ $F\in\PiM$ \ and \ $G\in\PhC$. 
\end{enumerate}
\end{defn}

\begin{defn}\label{dscoh}
Given a cohomological setting \ $\lra{\cC,\cM,\cV,\Phi,\Ap}$, \
take \ $\cE:=\sC$, \ with the resolution model category
structure defined by $\cM$. Then for any object $X$ and \Phal\ $G$ 
in $\cC$, the \emph{cohomology of $X$ with coefficients in $G$} 
is the total left derived functor of \ $\map_{\cC}(-,G)$, \ applied to
$X$. The \emph{$n$-th cohomology group of $X$ \ with coefficients in $G$} 
is the $\cM$-graded group:
$$
H^{n}(X;G):=[T^{n}\cu{\Ap M},(\LL\map_{\cC}(-,G)) X]_{M\in\cM}~.
$$
\end{defn}

\begin{defn}\label{dshom}
For \ $\lra{\cC,\cM,\cV,\Phi,\Ap}$ \ as above,
note that \ $\Ap M$ \ is a homotopy cogroup object in \ $\PhC$ \ for 
each \ $M\in\cM$, \ so we have a  resolution model category structure 
on \ $s\PhC$ \ determined by the set of models \ 
$\MP:=\{\Ap M\}_{M\in\cM}$. \ 
Define the \emph{homology} of $X$ to be the total left derived 
functor of \ $\Ap$ \ applied to $X$. \ The \emph{$n$-th homology group} of \
$X\in\cC$ is the $\cM$-graded group:
$$
H_{n}X~:=~\pin{\MP}{n}\LL\Ap X
$$
(cf.\ \S \ref{rsrmc}). \ (For this part of the definition we only require that \
$\PhC$ \ be enriched over itself via \ $\uHom$ \ -- \ we do not need
the symmetric monoidal structure.)

If \ $G\in\PhC$, \ we define the \emph{$n$-th homology group of $X$
  with coefficients in $G$} to be:
$$
H_{n}(X;G):=\pin{\MP}{n}(\LL(\Ap(-)\otimes G)(X)
$$
\end{defn}

\begin{example}\label{eshom}
The simplest example is when \ $\cC=\Gp$ \ (with \ $\cM=\{\ZZ\}$ \ as
on \S \ref{ermc}(i)), \ $\Phi=\fG$ \ (or $\fA$), and \ $\cV=\Set$, \ so \ 
$\PhC\cong\PhV\cong\Abgp$.

In this case \ $\PhC\cong\Abgp$, \ so the category \ $c\PhC$ \ of
cosimplicial $\Phi$-algebras in $\cC$ is equivalent to the category of
cochain complexes. Thus \ $K(\ZZ,n)$ \ (a cochain complex concentrated
in degree $n$) corepresents the $n$-th cohomology group of a cochain
complex \ ($n\in\NN$). \ This yields the usual cohomology groups of a
group $X$ with coefficients in an abelian group $G$ (as a trivial $X$-module).  

The functor \ $\Ap:\PiM\to\PhC$ \ is the abelianization \
$\Ab:\Gp\to\Abgp$, \ and the closed symmetric monoidal 
structure \ $\lra{\Abgp,\otimes,\ZZ,\Hom_{\Abgp}}$ \ yields the usual
homology of groups.
\end{example}

\begin{example}
\label{echsm}
Another simple example is provided by a symmetric monoidal category of
spectra, such as the symmetric spectra of \cite{HSSmiS}, or the
$S$-modules of \cite{EKMMayR}.  

In the latter version, for example, we take \ $\cE=\cM_{S}$, \ with
the symmetric monoidal smash product \ $\wedge_{S}$, and the internal function
complexes \ $F_{S}(-,-)\in\cV=\cE$ \ (cf.\ \cite[II, 1.6]{EKMMayR}). 
Since \ $\ho\cM_{S}$ \ is the usual stable homotopy category, it is
triangulated, with generator $S$. Thus we can take \ $\Phi=\ast$ \ to be
the trivial FP theory, any $S$-module $M$ yields a cohomology theory \
$F_{S}(-,M)$, \ and \ $\Ap:\cE\to\PhE$ \ is the identity. 
Similarly if \ $\cE=\cM_{R}$ \ for some $S$-algebra $R$.
\end{example}

\begin{remark}\label{rcohom}
These definitions may appear somewhat convoluted; they have been
set up to describe both the algebraic and (generalized) topological
theories in a uniform way, as appropriate derived functors.
Note that in general the total homology and cohomology functors, as
well as the homology and cohomology groups, take values in different
categories.
\end{remark}

%
%
\section{Theories and Abelianization}
\label{caft}

In this section we describe the necessary background for defining
(co)homology in a category \ $\cC=\TA$ \ of \Tal s. Most of it should
be familiar from the case \ $\cC=\Gp$, and the generalizations of Beck
and Quillen for algebras (see \cite{BecT,QuiC}); however, it seems
that the literature lacks a full description in this generality. We
start with the concept of (abelian) group objects, which are to play
the role of \Phal s in $\cC$. 

\begin{mysubsection}{Group objects}
\label{sgo}
In general, for a sketchable category \ $\cC=\TA$ \ we do not expect any
enrichment beyond \ $\cV=\Set$; \ so the natural choice for a
cohomological setting is \ $\Phi=\fA$.  

Recall that an \emph{(abelian) group object structure} on an object
$G$ in a category $\cC$ is a natural  (abelian) group structure on \ 
$\Hom\h_{\cC}(X,G)$ \ for all \ $X\in\cC$ \ -- \ in other words, a lifting 
of the functor \ $\Hom\h_{\cC}(-,G)$ \ from \ $\Set$ \ to \ $\Gp$ \
(or \ $\Abgp$); \ this is equivalent to a $\fG$-\ (resp., $\fA$-)
structure at $G$. 
\end{mysubsection}

\begin{remark}\label{rgo}
Note that if \ $\cC=\TA$ for some $\fG$-theory $\Theta$, any group
object structure on $G$ commutes with the underlying (graded) 
$\fG$-structure, so that the two necessarily agree and are commutative.  
In particular, in this case a \Tal\ can have \emph{at most one}
(necessarily abelian) group object structure. This is of course not
true for general $\cC$ (as is shown by the example of sets).
\end{remark}

\begin{mysubsection}{Abelianization of \Tal s}
\label{sata}
If $\Theta$ is any theory (sorted by $S$), the category of abelian
group objects in \ $\TA$ \ is sketched by the theory \ $\Tab:=\fA\Theta$ \ 
of \S \ref{sspat}. We call the $\fA$-localization \ $L_{\fA}:\TA\to\TAb$ \ 
the \emph{abelianization} functor for $\Theta$, and denote it
by \ $\Aa$. \ Note that \ $\Aa(\FT T)=\FTab T$.
\end{mysubsection}

\begin{examples}\label{egabel}
\begin{enumerate}\renewcommand{\labelenumi}{(\alph{enumi})}
\item When $\Theta$ is a $\fG$-theory, $\Tab:=\fG(\Theta)$, \ by
  Remark \ref{rgo}, and we can take \ $\Theta_{\fG}:=\Theta$ \ in \S
  \ref{sspat}, so \ $q:\Theta\epic\Tab$ \ is a quotient of theories, 
  and \ $q^{\ast}$ \ is simply the inclusion of the full subcategory
  of abelian \Tal s in \ $\TA$ \ (cf.\ \cite[\S 2.8]{BPescF}). \ Note that by 
  Remark \ref{renrich} we can then make \ $\Tab$ \ into a closed
  symmetric monoidal category. 
\item On the other hand, if \ $\Theta=\Td$, \ then \ $\Tab=\Theta_{\fA}$ \ 
  sketches $S$-graded abelian groups, \ $q^{\ast}:\TAb\to\TA$ \ is the
  forgetful functor \ $U:\gr_{S}\Abgp\to\gr_{S}\Set$, \ and its left
  adjoint \ $\Aa$ \ is the free graded abelian group functor\vsm.
\end{enumerate}
\end{examples}

\begin{mysubsect}{\Tal s over $X$}
\label{stalx}

We now show how the above discussion extends to the category \ $\TA/X$ \ 
of \Tal s over a fixed object $X$ (see \S \ref{rcomma}). First, we
need a:
\end{mysubsect}

\begin{defn}\label{dxalg}
If $\Theta$ is any theory and \ $X\in\TA$, \ then:
\begin{enumerate}
\renewcommand{\labelenumi}{(\alph{enumi})~}
\item An \emph{$X$-algebra} is an object $K$ in \ $\TA$ \ equipped
  with maps \ $\hat{f}:K(\vartheta)\times X(\vartheta)\to K(\vartheta')$ \ 
  for each \ $f:\vartheta\to\vartheta'$ \ in $\Theta$, satisfying: 
$$ 
\hat{g}(\hat{f}(k,x),X(f)(x))=\widehat{g\circ f}(k,x)
$$
for every \ $(k,x)\in K(\vartheta)\times X(\vartheta)$, \ and \ 
$g:\vartheta'\to\vartheta''$, \ with \ $\hat{f}(k,0)=K(f)(x)$.
\item The \emph{semi-direct product} of a \Tal\ $X$ by an $X$-algebra
  $K$ is the \Tal\ \ $K\rtimes X$ \ over $X$ given by: 
\begin{enumerate}
\renewcommand{\labelenumii}{(\roman{enumii})~}
\item $(K\rtimes X)(\vartheta):= K(\vartheta)\times X(\vartheta)$ \
  (as sets);
\item For each \ $f:\vartheta\to\vartheta'$ \ in $\Theta$, \ 
$(K\rtimes X)(f)(k,x):=(\hat{f}(k,x)),X(f)(x))$\vsm.
\end{enumerate}
\end{enumerate}
\end{defn}

If we want \ $K\rtimes X$ \ to be a group object in \ $\TA/X$, \ we
must require more. From now on, let $\Theta$ be a $\fG$-theory (sorted
by $S$), and $X$ a (fixed) \Tal. 

\begin{defn}\label{dxmod}
An \emph{$X$-module} is an $X$-algebra $K$ which is an abelian group
object in \ $\TA$, \ such that for each fixed \ $x\in X(\vartheta)$, \ 
each \ $\hat{f}(-,x):K(\vartheta)\to K(\vartheta')$ \ 
is additive (in the sense that it commutes with the given abelian
group structure). The category of $X$-modules will be denoted by \
$\RM{X}$ \ (see \cite[\S 3]{BecT}).
\end{defn}

\begin{remark}\label{rxmod}
In this case the underlying $S$-graded group \ $VK$ \ is an $S$-graded \ 
$VX$-module in the traditional sense (a module over the graded  group
ring \ $\ZZ[VX]$), \ and the group operation at each \ $\theta\in\Theta$ \ 
is given by \ $\mt((k,x),(\ell,y)~=~(k+x\cdot\ell,xy)$, \ 
as usual.
\end{remark}

\begin{defn}\label{dder}
Assume that \ $p:Y\to X$ \ is a map of \Tal s, and $K$ is an
$X$-module. A function \ $\xi:Y\to K$ \ (preserving the products of
$\Theta$) will be called a \emph{derivation with respect to $p$} if \ 
$\xi(Y(f)(y))=\hat{f}(\xi(y),p(y))$ \ for any \ 
$f:\vartheta\to\vartheta'$ \  in $\Theta$. \ The set of all such will
be denoted by \ $\Der\h_{p}(Y,K)$. \ 
In particular, a derivation with respect to \ $\Id\h_{X}$ \ will be
called simply a \emph{derivation}, and \ $\Der(X,K):=\Der\h_{\Id}(X,K)$. 
\end{defn}

\begin{remark}\label{rder}
Note that this holds in particular for \
$f=\mt:\theta\times\theta\to\theta$, \ so that by Remark \ref{rxmod}:
$$
\xi(\mt(y_{1},y_{2})))~=~
\hm((\xi(y_{1}),p(y_{1})),(\xi(y_{2}),p(y_{2}))~=~
\xi(y_{1})+p(y_{1})\cdot\xi(y_{2})~.
$$
Thus $\xi$ is a derivation (crossed homomomorphism) with respect
to the $\fG$-structure.

Furthermore, \ $\Der_{p}(Y,K)$ \ is an abelian group (under the
addition of $K$), and any map of $X$-modules \ $\alpha:K\to L$ \
induces a homomorphism \ $\alpha_{\ast}:\Der_{p}(Y,K)\to\Der_{p}(Y,L)$.
\end{remark}

The following results do not appear in this form in the
literature, but their proofs are straightforward generalizations of
the corresponding (classical) results for groups (see, e.g., 
\cite[\S 3-4]{BecT}  and \cite[\S 11.1]{RobiG}).

\begin{prop}\label{pmodgp}
Any group object structure on \ $p:Y\to X$ \ in \ $\TA/X$ \ is necessarily
abelian. Moreover, \ $K:=\Ker(p)$ \ is an $X$-module, with \
$Y\cong K\rtimes X$, \ and for some derivation \ $\xi:X\to K$, \ the group
operation map \ $\mu:Y\times_{X}Y\to Y$ \ is given  (under the
identification \ $UY=UK\times UX$) \ by \
$\mu(k,k',x)=(k+k'+\xi(x),x)$, \ the zero map by \ $(k,x)\mapsto
(-\xi(x),x)$, \ and the inverse by \ $(k,x)\mapsto(-k-2\xi(x),x)$. 

Conversely, for any  $X$-module $K$ and derivation \ $\xi:X\to K$, \
the above formulas make \ $K\rtimes X$ \ into an abelian group object
over $X$.
\end{prop}

\begin{cor}\label{cabgpobj}
There is an equivalence of categories \ $\ell^{\ast}\hy{\fG}{(\TA/X)}\to\ATAX$, \ 
induced by the quotient map \ $\ell:\fG\hra\fA$. 
\end{cor}

\begin{lemma}\label{lhomom}
Any homomorphism \ $\phi:K\rtimes X\to L\rtimes X$ \ between group
objects over $X$ (with group operations determined by \ 
$\sigma\in\Der(X,K)$ \  and \ $\tau\in\Der(X,L)$, \ respectively) is 
of the form \ $\phi(k,x)=(\alpha(k)+\xi(x),x)$, \ where \ 
$\alpha:K\to L$ \ is a homomorphism of $X$-modules and \ 
$\xi:=\alpha\circ\sigma-\tau$.
\end{lemma}

In particular, any two group object structures over $X$ on the
semi-direct product \ $K\rtimes X$ \ are canonically isomorphic, so
we deduce:

\begin{prop}\label{pcorresp}
The functor \ $\lambda:\RM{X}\to\ATAX$, \ defined \ 
$\lambda(K):=K\rtimes X$ \ (with the group operation map determined by
the zero derivation), is an equivalence of categories, with inverse \
$\kappa:\ATAX\to\RM{X}$ \ which assigns to an abelian group
object \ $p:Y\to X$ \ the $X$-module \ $\Ker(p)$.
\end{prop}

\begin{remark}\label{rsdp}
Since the forgetful functor \ $U=\UT:\TA\to\TdA$ \ is faithful, for any
\Tal\ $Y$ and semi-direct product \ $K\rtimes X\in\TA$ \ we have:
%
\begin{myeq}\label{efour}
\begin{split}
&\Hom\h_{\TA}(Y,K\rtimes X)~\stackrel{U}{\hra}~
\Hom\h_{\TdA}(UY,U(K\rtimes X))~=\\
&\Hom\h_{\TdA}(UY,UK\times UX)~=~
\Hom\h_{\TdA}(UY,UK)\times\Hom\h_{\TdA}(UY,UX)~.
\end{split}
\end{myeq}
Thus given \ $p:Y\to X$, \  we can write any map \ $\phi:Y\to K\rtimes X$ \ 
over $X$ \ in the form \ $\phi(y)=(\alpha(y),p(y))$, \ and the
requirement that $\phi$ be a map of \Tal s means that \ 
$\alpha:\FT T\to K$ \ is a derivation with respect to $p$
(\S \ref{dder}), so in fact: 
%
\begin{myeq}\label{eseventeen}
\Hom\h_{\TA/X}(Y,K\rtimes X)~\cong~\Der_{p}(Y,K)
\end{myeq}
\noindent as abelian group (once we choose a fixed group structure
on \ $K\rtimes X$)\vsm.

Three special cases should be noted:

\begin{enumerate}\renewcommand{\labelenumi}{(\alph{enumi})}
\item For \ $p=\Id:X\to X$, \ we see that \ $\Der(X,K)$ \ is the space
  of sections for \ $K\rtimes X$, \ as usual.
\item If \ $Y=L\rtimes X$ \ for some \ $L\in\RM{X}$, \ then by
  Proposition \ref{pcorresp}:
$$
\Hom\h_{\RM{X}}(L,K)~\stackrel{\lambda}{=}~
\Hom\h_{\TA/X}(L\rtimes X,K\rtimes X)~=~\Der_{p}(L\rtimes X,K)~.
$$
On the other hand, by Lemma \ref{lhomom} any map of $X$-modules \ 
$\alpha:L\to K$ \ induces a homomorphism of group objects \ 
$\phi=\lambda(\alpha):L\rtimes X\to K\rtimes X$ \ (where we use the
zero derivation to define the group structures on the semi-direct
products). Thus in fact:
%
\begin{myeq}\label{eeighteen}
\Hom\h_{\ATAX}(L\rtimes X,K\rtimes X)~=~\Der_{\pi_{2}}(L\rtimes X,K)
\end{myeq}
\noindent as abelian groups.
\item If \ $Y=\FT T$ \ is free, then by adjointness we actually have
  equalities of sets:
$$
\Hom\h_{\TA}(\FT T,K\rtimes X)=\Hom\h_{\TdA}(T,UK)\times\Hom\h_{\TdA}(T,UX)
$$
in \ \eqref{efour}, \ so for \ $p:\FT T\to X$ \ in \ $\fFT/X$, \ 
we have: 
%
\begin{myeq}\label{enineteen}
\Hom\h_{\TA/X}(\FT T,K\rtimes X)~\cong~\Hom\h_{\TdA}(T,UK)~\cong~
\Hom\h_{\TA}(\FT T,WK)~,
\end{myeq}
\noindent where \ $W:\RM{X}\to\TA$ \ is the forgetful functor.
In particular:
$$
\Der_{p}(\FT T,K)~\cong~\Hom\h_{\TA}(\FT T,WK)
$$
\noindent as sets (though this identification is not natural in the
full subcategory \ $\fFT$ \ in \ $\TA$). 
\end{enumerate}
\end{remark}

\begin{mysubsection}{Abelianization over a \Tal}
\label{satax}
Recall from \S \ref{rcomma} that for a fixed \Tal\ $X$, \ $\TA/X$ \ 
can be sketched by \ $\ToX$ \ (sorted by \ $\UT X$). \  Similarly, \
$\ATAX$ \ can be sketched by \ $\AToX$, \ obtained from \  $\ToX$ \ as
in \S \ref{sspat} by adding:

\begin{enumerate}\renewcommand{\labelenumi}{(\alph{enumi})}
\item a section \ -- \ i.e., constants in each \ $\phi_{x}$ \ 
(in the notation of \S \ref{rcomma});
\item group structure maps \ $\mu:\phi_{x}\times\phi_{x}\to\phi_{x}$ \
  and \ $\rho:\phi_{x}\to\phi_{x}$, 
\end{enumerate}
\noindent satisfying the obvious identitites. Again the map of theories \ 
$i:\ToX\hra\AToX$ \ induces the forgetful functor \ 
$i^{\ast}:\ATAX\to\TA/X$, \ with an adjoint \ $\Aox:\TA/X\to\ATAX$ \  
called the \emph{abelianization} of \ $\TA/X$. \ This is needed in
order to define homology for \Tal s (see \S \ref{shomt} below).

Note that the category \ $\RM{X}$ \ can also be sketched by an
$\fA$-theory \ $\TX$, \ obtained from \ $\Tab$ \ (\S \ref{sata}) by
adding operations \ $x\cdot(-):\theta\to\theta$ \ for each \ $x\in \UT
X$, \ satisfying  the obvious identitites. \ The inclusion \ $j:\Tab\hra\TX$ \ 
induces the forgetful functor \ $j^{\ast}:\TXA\to\TAb$. \ 
If we define \ $\kappa:\TA/X\to\TA$ \ as in Proposition \ref{pcorresp}, 
we obtain the commutative outer diagram:
$$
\xymatrix@R=25pt{
\TAXa \ar[rr]^{i^{\ast}} \ar@<-2.5ex>[d]_{\kappa} & &
\TA/X \ar[d]^{\kappa}\ar@{-->}[lld]_{\hAox}\\
\RM{X}=\TXA\ar[u]_{\lambda}^{\cong} \ar[r]_>>>>{j^{\ast}} & 
\TAb\ar[r]_{q^{\ast}} & \TA
}
$$
in which the horizontal arrows are forgetful functors (and \
$q^{\ast}$, \ $i^{\ast}$ \ have adjoints \ $\Aa$, \ $\Aox$, \
respectively, \ with \ $\hAox:=\kappa\circ\Aox:\TA/X\to\RM{X}$).

Note that by \ \eqref{enineteen}, \ the abelianization functor \ 
$\hAox$ \ takes any free \Tal\ \ $p:\FT T\to X$ \ over $X$ to the
corresponding free $X$-module \ $\FTx T\in\TXA=\RM{X}$. \ Moreover, for any \ 
$\varphi\in\Der_{p}(\FT T,K)$ \ (determined by \
$\varphi(t_{i})=k_{i}\in K$ \ for \ $t_{i}\in T$), \ the corresponding \ 
$\hat{\varphi}\in\Hom\h_{\RM{X}}(\FTx T,K)$ \ is also determined by 
requiring that \ $\hat{\varphi}(t_{i})=k_{i}$. \ 
Now assume given a map \ $\psi:\FT T'\to\FT T$ \ in \ $\fFT/X$, \ 
determined by the condition that, for each \ $t'\in T'$,\ 
$\psi(t')=f'_{\ast}(t_{i_{1}},\dotsc,t_{i_{n}})$ \ for some $f'$ in
$\Theta$. Then:
$$
(\psi^{\ast}\varphi)(t')~=~
\hat{f}'((t_{i_{1}},\dotsc,t_{i_{n}}),(p(t_{i_{1}}),\dotsc,p(t_{i_{n}})))\in K~.
$$
\end{mysubsection}

\begin{remark}\label{rtrivmod}
Evidently, the discussion of abelian group objects and abelianization
over a \Tal\ $X$ extends the absolute case of \S \ref{sgo}\textit{ff.}, 
taking \ $X=0$. 

More generally, $K$ will be called a \emph{trivial}
$X$-module if \ $\hat{f}(k,x)=f(k)$ \ for every \ $f\in\Theta$ 
(\S \ref{dxalg}) \ -- \ so that $K$ is simply an abelian \Tal, \ 
$K\rtimes X$ \ is the product in \ $\TA$, \ and a derivation into $K$
is just a map of \Tal s.  
\end{remark}

%
%
\section{(Co)homology of \Tal s}
\label{cch}

Andr\'{e} (in \cite{AndrM}) and Quillen (in \cite[II, \S 5]{QuiH} and 
\cite[\S 2]{QuiC}) defined homology and cohomology groups in
categories of universal algebras. Quillen also showed how this
generalized the earlier definition of triple cohomology (see 
\cite[\S 2]{BecT}). We now indicate briefly how this definition fits
into the setup of \S \ref{scohf}.    

\begin{mysubsection}{Cohomology of \Tal s}
\label{scoht}
Let $\Theta$ be a $\fG$-theory, and \ $\cC:=\TA$ \ (or \ $\TA/X$ \ for a
fixed \Tal \ $X$), \ with the resolution model category structure on \ 
$\sC$ \ described in \S \ref{ermc}(ii) (or \S \ref{rcomma}).

As in Example \ref{eshom}, here \ $\cV=\Set$, \ so we must take \ $\Phi=\fA$ \  
(or equivalently, by Corollary \ref{cabgpobj}: \ $\Phi=\fG$), \ 
since cosimplicial \emph{sets} do not have any useful model category
structure (see however \cite{BouC}). Thus if $G$ is an abelian group
object in $\cC$, and \ $\Vd\to Y$ \ is a free simplicial resolution
(cofibrant replacement in \ $\sC$), then the cosimplicial abelian
group \ $\Wu:=\Hom\h_{\cC}(\Vd,G)$ \ corresponds under the Dold-Kan
equivalence (cf.\ \cite[\S 3]{DPupH} and \cite[9.4]{WeibHA}) to a
cochain complex \ $\Wus$, \ and the category \ $\cCh{\ZZ}$ \ of
non-negatively graded cochain complexes of abelian groups embeds in
the category \ $\Ch{\ZZ}$ \ of unbounded (co)chain complexes, which is
a stable model category (cf.\ \cite[Ch.\ 7]{HovM}. \ Suspensions of \
$g:=K(\ZZ,0)$ \ detect homology in \ $\cCh{\ZZ}$ \ (or \ $\Ch{\ZZ}$),
\ so \ $c\Abgp\cong\cCh{\ZZ}$ \ is semi-triangulated in the sense of
\S \ref{dqtc}, \ and in fact \ $[T^{i}g,\Wus]=H^{i}(Y;G)$ \ is the
$i$-th  Andr\'{e}-Quillen cohomology group of $Y$. 

Remark \ref{rsdp} shows that these can be thought of as usual as the
derived functors of \ $\Der(-,G)$, \ in the case \ $\cC=\TA/X$, \ and
as \ $\Ext^{i}(Y,G)$ \ in the case \ $\cC=\TA$ \ (\S \ref{rtrivmod}).
This identification has been the basis for a number of definitions of
cohomology in various topological settings \ - \ see, e.g., \cite{MShniCD},
and the survey in \cite{BRichtC}.
\end{mysubsection}

\begin{mysubsection}{Homology of \Tal s}
\label{shomt}
In this situation one can define the homology of a \Tal\ $Y$ as the
total left derived functor of abelianization \ $\Aa:\TA\to\TAb$ \ 
(\S \ref{sata}), \ which takes values in the category \ $s\TAb$ \ of
simplicial \ $\Tab$-algebras (as usual, we only need to evaluate \ 
$\Aa$ \ on \ $\fFT$, \ so \ $\LL\Aa$ \ actually takes values in \ 
$s\fF_{\Tab}$). \ Since \ $\TAb$ \ is an abelian category (with enough 
projectives, namely: \ $\fF_{\Tab}$), \ $s\TAb$ \ is equivalent
to the stable model category \ $\ch(\Tab)$ \ of chain complexes over \
$\Tab$, \ and the homology groups \ $[T^{i}K(\FTab s,0),\Aa\Vd]=H_{i}Y$ \ 
(for $s$ an $S$-graded singleton) are themselves \ $\Tab$-algebras. 

The same holds for \ $Y\in\TA/X$: \ using \S \ref{satax}, we may
define \ $H_{i}(Y/X)$ \ as the $i$-th derived functor of \
$\Ax:\fFT/X\to\TAXa$, \ taking values in \ $\TAXa$ \ -- \ or
equivalently (Proposition \ref{pcorresp}) in $X$-modules. For groups, \ 
$H_{\ast}(G/G)$ \ is the homology of $G$ with coefficients in \
$\ZZ[G]$. \ For a pointed connected space $X$ with \ $G=\pi_{1}(X,x)$, \ 
$H_{\ast}(X/BG)$ \ is the homology of $X$ with coefficients in the
local system \ $\ZZ[G]$. 
\end{mysubsection}

\begin{defn}\label{dhomcoef}
To define homology of \ $Y\to X$ \ with coefficients in an arbitrary
$X$-module $G$, \ we need a monoidal structure on \
$\RM{X}\cong\TAXa$, \ induced via the adjoint pair
$$   
\TXA\adj{\FTx}{\UTx}\TdA
$$
from the usual monoidal structure \ $(\TdA,\times)$ \ of Cartesian products of
graded sets.

More precisely, define \ $\otimes:\fF_{\TX}\times\fF_{\TX}\to\fF_{\TX}$ \ 
by \ $\FTx T\otimes\FTx S:=\FTx(T\times S)$. \ The $0$-th derived
functor in the second variable defines \ $\FTx T\otimes G$ \ for any
$\TX$-algebra ($X$-module) $G$; and the $n$-th left derived functor of \
$\Ax(-)\otimes G$ \ (in the first variable) is by definition \
$H_{n}(Y/X;G)$. 
\end{defn}

\begin{example}\label{eghomology}
When \ $\Theta=\fG$, \ a free simplicial resolution \ $\Vd$ \ of a
group $G$ in \ $s\Gp$ \ is actually a cofibrant model for the
classifying space \ $BG$ \ (in \ $\Sa$). \ Applying the functor \
$\hAox$ \ of \S \ref{satax} to \ $\Vd$ \ dimensionwise yields a model
for the chains on the universal contractible $G$-space \ $EG$ \ 
(since conversely, taking the free $\ZZ$-module on the bar
construction model for \ $EG$ \ and dividing out by the free
$G$-action yields \ $\ZZ BG$, \ so \ $\ZZ EG\simeq \ZZ[G]\Vd$). \ 
Taking homotopy groups of \ $\ZZ[G]\Vd$ \ is the same as taking the
homology of the chain complex corresponding to \ $\ZZ EG$, \ which
is just \ $H_{\ast}(G;\ZZ[G])$.
\end{example}

\begin{remark}\label{rtops}
Note that the previous discussion actually defines homology and
cohomology for any simplicial \Tal\ \ $\Yd$, \ not only for the
constant ones. Moreover, if \ $\Theta=\fG$, \ the adjoint pairs of functors:
%
\begin{myeq}\label{ethree}
\Ta\adj{|-|}{S}\Sa\adj{G}{\bar{W}}\cG=s\Gp
\end{myeq}
\noindent induce equivalences of the homotopy categories of pointed
connected topological spaces, reduced simplicial sets, and simplicial
groups. Here \ $|-|$ \ is the geometric realization functor, $S$ is the
singular set functor, $\bar{W}$ is the Eilenberg-Mac~Lane
classifying space functor, and $G$ is Kan's loop functor 
(cf.\ \cite[\S 26.3]{MayS} and \cite[I,4 \& II,3]{QuiH}). 
Thus Quillen's approach provides an algebraic description
of ordinary homology and cohomology of spaces (with local
coefficients). Note, however, the shift in indexing: in particular, we
lose \ $H^{0}$, \ since we can deal only with connected spaces from
this point of view.  

There is also an algebraic model for not-necessarily-connected spaces
due to Dwyer and Kan, using simplicial groupoids (see \cite[V, \S 5]{GJarS}), 
and Quillen's approach, as well as much of the discussion here,
carries over to that setting (compare \cite{DuskN}).  However, in
order to avoid further complicating the description, we restrict
attention here to simplicial groups. 
\end{remark}

\begin{mysubsection}{Diagrams of \Tal s}
\label{sdiagt}
If $D$ is a small category and $\Theta$ is a $\fG$-theory, there
is a model category structure on the functor category \ $s\TA^{D}$,
and the objectwise descriptions of abelian group objects and
abelianization (for each \ $d\in D$) \ provide definitions of
(co)homology for diagrams of \Tal s, too (see \cite[\S 4]{BJTurR} for
the details). 

Moreover, even for \ $\cC=\TA$ \ or \ $\TA/X$, \ we can allow our 
coefficients to be \emph{diagrams}  \ $G:D\to\ATA$ \ of abelian group
objects (or $X$-modules). This enables us to treat a map such as \ 
$\ZZ\epic\ZZ/p$ \ (reduction mod $p$), say, as the coefficients for a
cohomology theory (rather than a natural transformation). In
particular, we can apply any general machinery, such as universal
coefficient theorems, to \ $H^{\ast}(-;G)$, \ too\vsn.
\end{mysubsection}

\begin{mysubsect}{Spherical model categories}
\label{ssmc}

When \ $\cC=\TA$ \ for some $\fG$-theory $\Theta$, the resolution model
category \ $\sC$ \ (and the models \ $\cM=\fFTp$ \ - \ cf.\ 
\S \ref{ermc}(ii)) \ will have additional useful structure which is
familar to us from topological spaces\vsm : 

\begin{enumerate}
\renewcommand{\labelenumi}{\arabic{enumi}.~}
\item For any \ $n\geq 1$, \ $\picM{n}(-)$ \ is naturally an
  abelian group object over \ $\picM{0}(-)$. 
\item  Each \ $\Vd\in \sC$ \ has a functorial \emph{Postnikov tower} 
of fibrations: 
$$
\dotsc \to P_{n}\Vd\xra{p\q{n}}P_{n-1}\Vd\xra{p\q{n-1}}\dots\to P_{0}\Vd~,
$$
as well as a weak equivalence \ 
$r:\Vd\to P_{\infty}\Vd:=\lim_{n}P_{n}\Vd$ \ 
and fibrations \ $P_{\infty}\Vd\xra{r\q{n}}P_{n}\Vd$ \ such that \
$r\q{n-1}=p\q{n}\circ r\q{n}$ \ for all $n$, \ and \ 
$(r\q{n}\circ r)_{\#}:\picM{k}\Vd\to\picM{k}P_{n}\Vd$ \ 
is an isomorphism for \ $k\leq n$, \ and zero for \ $k>n$. \ 
\item For every \Tal\ $X$, there is a \emph{classifying object} \ 
$BX$ \ with \ $BX\simeq P_{0}BX$ \ and \ $\picM{0}BX\cong X$, \
unique up to homotopy. 
\item Given a \Tal\ $X$ and an $X$-module $G$, there is an 
\emph{extended $G$-Eilenberg-Mac~Lane object} \ $E=E^{X }(G,n)$ \ 
in \ $\sC/X$ \ for each \ $n\geq 1$, \ unique up to homotopy, 
equipped with a section $s$ for \ $p\q{0}:E\to P_{0}E\simeq BX$, \ 
such that \ $\kappa\picM{n}E\cong G$ \ as $X$-modules; \ and \ 
$\picM{k}E=0$ \ for \ $k\neq 0,n$. \ If $G$ is a trivial
$X$-module (\S \ref{rtrivmod}), we write simply \ $E(G,n)$.
\end{enumerate}

Any resolution model category with this additional structure (as well
as functorial $k$-invariants) is called a \emph{spherical model
category}. See \cite[\S 1-2]{BlaC} for the details.
\end{mysubsect}

\begin{remark}\label{rsphx}
The homotopy groups \ $\piM{n}$ \ in the resolution model category \
$s\TA$ \ are corepresented by \ $S^{n}\otimes \FT(s)$ \ for \ 
$M=\FT(s)\in\fFTp$, \ $s\in S\subseteq\Theta$ \ (cf.\ \S \ref{rsrmc}).
Thus by adjointness for any \ $\Vd\in s\TA$ \ we have:
$$
\piM{n}\Vd~=~[S^{n}\otimes \FT(s),\Vd]_{\ast}~=~[S^{n},(\UT \Vd)_{s}]~=~
\pi_{n}(\UT \Vd)_{s}~,
$$
\noindent so that the group \ $\piM{n}X$ \ (induced by
the homotopy cogroup structure of \ $S^{n}$) \ is the usual $n$-th 
simplicial homotopy group of the graded simplicial group \ $\UT\Vd$ \ 
in the appropriate degree.

This works also in \ $s\TA/X$: \ more precisely, \ $\picM{n}\Vd$ \ as
defined above is an abelian group object in \ $\TA/\picM{0}\Vd$, \ and
applying $\kappa$ of Proposition \ref{pcorresp} yields a \
$\picM{0}\Vd$-module, whose underlying $S$-graded set is \
$\pi_{n}\UT\Vd$ \ (see \cite[\S 4.14]{BPescF}).
\end{remark}

\begin{mysubsection}{Cohomology in \ $s\TA$}
\label{scohsta}
It may appear more natural to take as a representing object an abelian
group object in the model category \ $s\TA$ \ itself. In most cases
this will yield no new cohomology groups, but it will enable us
to define, and in some cases compute, the primary cohomology
\emph{operations} \ -- \ as we do for topological spaces 
(see, e.g., \cite{PridP}).

The obvious examples are those of the form \ $E(G,n)$ \ as above (or \
$E^{X}(G,n)$ \ in \ $s\TA/BX $, \ if we want local coefficients). In most 
cases of interest \ -- \ including \ $\Ta$, \ $\Sa$, \ $\cG=s\Gp$ \ -- \ 
the only objects  in \ $\hy{\fA}{s\TA}$ \ are products of the above.
Furthermore, since \ $E(-,n):\ATA\to s\TA$ \ is a functor, \ we can
define an Eilenberg-Mac~Lane diagram \ $E(G,n)$ \ for any
diagram \ $G:D\to\ATA$ \ as in \S \ref{sdiagt}.

Thus for any cofibrant \ $\Wd$ \ in \ $s\TA$ \ and coefficients \ 
$\cM\in\ATA^{D}$, \ for each \ $n\geq 1$ \ we define the $n$-th
\emph{cohomology group} of \ $\Wd$ \ with coefficients in $G$,
denoted by \  $H^{n}(\Wd;G)$, \  to be the set of components of \
$\map(\Wd,E(G,n))$ \ (which is a $D$-diagram of simplicial abelian
groups, so the components constitute a $D$-diagram of abelian groups). 

Again, there is also a local version, for $G$ in \ $\TA/X$ \ or \ 
$\cM:D\to\ATA/X$, \ yielding:
$$
H^{n}(\Wd/X;G):=\pi_{0}\map\h_{s\TA/X}(\Wd,E^{X}(G,n))
\hsp \text{for each} \ n\geq 1.
$$
\end{mysubsection}

\begin{prop}\label{pone}
If $\Theta$ a $\fG$-theory, $X$ is a \Tal, and $G$ is in \ $\ATAX$, \ 
then cohomology with coefficients in $G$ as defined in \S \ref{scoht}
is naturally isomorphic to that defined in \S \ref{scohsta}.
\end{prop}

Compare \cite[\S 3]{DuskS}.

\begin{proof}
Let $K$ be the $X$-module corresponding to \ $G=K\rtimes X$, \ so \ 
$\Ed:=E^{X}(G,n)$ \ is of the form \ $E(K,n)\rtimes X$, \ where \  
$E(K,n)$ \ is obtained from the analogous chain complex (over \
$\RM{X}$) \ by the Dold-Kan equivalence (cf.\ \cite[p.~95]{MayS}). Thus:
%
\begin{myeq}\label{esix}
E_{i}~=~\begin{cases}
             X ~                               & \text{for}\hsm 0\leq i<n\\
             K\rtimes X                            &\hs i=n\\
             (\bigoplus_{j=0}^{n}~s_{j}K)~\rtimes X ~&\hs i=n+1\\
             M_{i}\Ed~                                 &\hs i\geq n+2~,
            \end{cases}
\end{myeq}
\noindent (where \ $M_{i}\Ed$ \ is the $i$-th \emph{matching object} \ -- \ 
see \cite[X, \S 4.5]{BKaH} or \cite[\S 2.1]{BJTurR}), with the differential:
%
\begin{myeq}\label{eseven}
\partial_{n+1}(x,\lambda))~:=~(\sum_{i=0}^{n+1}~d_{i}x,\lambda)\hs 
\text{for every}\hsm  (x,\lambda)\in E_{n+1}~.
\end{myeq}

Let \ $\Wd$ \ be a free simplicial object in \ $\sC$, \ with \ 
$\varepsilon:W_{0}\to X$ \ inducing \ $\pi_{0}\Wd\cong X$ \ 
(for example, \ $\Wd$ \ could be a resolution of $X$).
From \eqref{esix} and \eqref{eseven} we see that \
$\Hom\h_{\sC/X }(\Wd,\Ed)$ \ is naturally isomorphic to the 
subgroup of \ $\Hom\h_{\cC/X}(W_{n},K\rtimes X)$ \
consisting of maps \ $f:W_{n}\to K\rtimes X$ \ (over $X$)
for which \ $f\circ d_{i}$ \ is the projection to $X$ \ (the
zero of $\Hom(\cC/X )(W_{n+1},K\rtimes X)$ \ for each \ 
$0\leq i\leq n+1$. \ Here \ $W_{n}$ \ maps to $X$ by \ 
$\varepsilon\circ d_{0}\circ\dotsb\circ d_{0}$. \ 

Again by the Dold-Kan equivalence, there is a path object \ $\Ed^{I}$ \
for \ $\Ed$ \ in \ $s\TA/X$ \ with
%
\begin{myeq}\label{eeight}
E^{I}_{i}~=~\begin{cases}
             X ~                               &\text{for}\hsm0\leq i<n-1\\
             K\rtimes X                         &\hs i=n-1\\
             (K\oplus K\oplus\bigoplus_{j=0}^{n-1}~s_{j}K)~\rtimes X~  
                                                    &\hs i=n\\
             M_{i}\Ed~                              &\hs i\geq n+1~,
            \end{cases}
\end{myeq}
\noindent with \ $d_{0}$ \ the identity on the first copy of \
$K\rtimes X$ \ in \ $E^{I}_{n}$, \ and minus the identity on the
second copy. There are two obvious projections \ $p_{0},p_{1}:\Ed^{I}\to\Ed$, \ 
and a homotopy between two maps \ $f_{0},f_{1}:\Wd\to\Ed$ \ over $X$ 
is a map \ $F:\Wd\to\Ed^{I}$ \ with \ $p_{i}\circ F=f_{i}$ \ $(i=0,1$), \ 
which in turn corresponds to a map \ $F':W_{n-1}\to K\rtimes X$ \ 
over $X$ for which \ 
$F'\circ d_{0}$ \ represents \ $f_{0},f_{1}$ \ respectively on the two
copies of \ $K\rtimes X$. \ 

Thus we see that \ $H^{n}(\Wd/X,M):=[\Wd,\Ed]_{\sC/X}$ \ 
is canonically isomorphic to the $n$-th cohomotopy group of the
cosimplicial abelian group \ $\Hom(\cC/X)(\Wd,K\rtimes X)$, \ 
as claimed.
\end{proof}

\begin{mysubsect}{Cohomology of operads and their algebras}
\label{scoa}

As noted in \S \ref{egtheory}(b), our definition of sketchable
categories covers both the category of operads, \ $\Alg{\cO}$, \ and
that of algebras over a given operad $\cP$. 

Of course, $\cO$ is not a $\fG$-theory; however, essentially all known
applications are to operads of (connected) topological spaces or of
chain complexes (see \cite{MSStasO}). In the first case, we can use \
\eqref{ethree} \ to replace \ $\Ta$ \ by $\cG$, so that in both cases
we may assume, without loss of generality, that our operad takes value
in \ $s\TA$ \ for some $\fG$-theory $\Theta$.  
Note that the category of $\cO$-algebras in \ $s\TA$ \ is equivalent
to \ $s\TtA$, \ where \ $\tilde{\Theta}=\cO\times\Theta$ \ (product of
FP-sketches) is now an $\fG$-theory (see \S \ref{sspat}). Thus the definition  
of \S \ref{scohsta} (applied to $\tilde{\Theta}$) is valid for operads
of spaces or chain complexes.

The same applies to algebras over a fixed operad $\cP$ taking values
in \ $\Ta$ \ or \ $\Ch{k}$ \ for some field $k$ (see \cite[\S 2]{MayG}), 
as well as to the cohomology of a $k$-linear category (that is,
algebras over a $k$-linear PROP) considered in \cite{MarkCC}.

We should observe, however, that the various cohomology theories
constructed  \ -- \ in the context of deformation theory \ -- \ in
\cite{MarkCC}, in \cite{MShniDA} (for Drinfel'd algebras), in
\cite{GSchaB} (for bialgebras), and so on, are defined in terms of a
specific differential graded resolution. To show that these agree with
our general definition requires a generalization of Quillen's
equivalence between simplicial and differential graded Lie algebras
over $\QQ$ (see \cite[I, \S 4]{QuiR}, and compare \cite[\S 3]{DPupH}). 
One can expect such an equivalence only for suitable $k$-linear
categories over a field $k$ of characteristic $0$.  
\end{mysubsect}

\begin{remark}\label{rjpir}
We should point out that a different definition of (co)homology for
\Tal s, based on the Baues-Wirsching and Hochschild-Mitchell cohomologies
of categories  (cf.~\cite{BWirC,BMitR}), \ is given by Jibladze and
Pirashvili in \cite{JPirCA}. See \cite[Theorem~6.7]{SchweSA} for
an equivalent formulation in terms of the topological Hochschild
(co)homology of suitable ring spectra.
\end{remark}

\begin{mysubsection}{Cohomology of sheaves}
\label{scohsh}
We have assumed so far that $\Theta$ was a $\fG$-theory. This is
necessary for the approach described here at two points: in order 
to identify the (abelian) group objects in \ $\TA$ \ (see Section
\ref{caft}), and to define the model category structure on \ $s\TA$ \
(see \S \ref{ermc}(ii)). This is a resolution model category (induced
by the adjoint pair \ $(\FT,\UT)$ \ of \S \ref{sfta}) only with some
such additional assumption (cf.\ \cite{BlaN}): otherwise the free \Tal s 
are not necessarily  cogroup objects. 

One obvious example where this fails is the category of sets, where we
apparently have no meaningful concept of cohomology. A more
interesting case is the category of sheaves on a topological 
space $X$, sketched by \ $\TU$ \ (see \S \ref{egtheory}). Note
that there is no free/forgetful adjoint pair between \ 
$\Alg{\Theta^{\delta}_{\cU}}$ \ and \ $\TUA$ \ or \ 
$\Tab=\hy{\fA}{\TU}\cong\hy{\TU}{\Abgp}$, \ since sheaves of abelian
groups rarely have \emph{any} projectives (e.g., \ $\ZZ C_{U}$ \ 
in \S \ref{egtheory}\,(c) is not generally a sheaf).   
However, they do have enough injectives, so if we replace
left derived functors by right derived functors in \S \ref{scohf},
with \ $\cE=\TUA$, \ $\cV=\Set$, \ and \ $\Phi=\fA$, \ we may define \ 
$H^{n}(X;\cF)$, \ for any \ $\cF\in\PhE$, \ to be the right derived
functors of \ $\Hom\h_{\,E}(C_{X},-)$, \ applied to $\cF$. \ 
This also explains why our definition of homology does not make sense
for sheaves. 
\end{mysubsection}

%
%
\section{Generalized cohomology}
\label{cgencoh}

For simplicial \Tal s over a $\fG$-theory $\Theta$  \ -- \ and thus
for simplicial sets or topological spaces \ -- \  the only strict
abelian group objects are generalized Eilenberg-Mac~Lane objects 
(cf.\ \cite[19.6]{MoorH2}). Of course, in any model category $\cDD$,
any abelian group object $G$ in \ $\ho\cDD$ \ defines a functor \
$[-,G]:\ho\cDD\to\Abgp$; but such functors do not usually satisfy the axioms
of a cohomology theory. From our point of view, this is because 
the structure maps on the higher products \ $G^{k}$ \ ($k\geq 3$) \
which are needed to make $G$ an $\fG$- or $\fA$-algebra in $\cDD$ are
not uniquely defined. 

One way to deal with this problem would be to require that $G$ have an \
$E_{\infty}$-operad acting on it (cf.\ \cite[\S 14]{MayG}).  
If \ $\cDD=\Ta$ \ (or \ $\Sa$), by a result of Boardman and Vogt,
under mild topological restrictions any \ $E_{\infty}$ \ $H$-space is
homotopy equivalent to a strict abelian monoid in $\cDD$ (cf.\ 
\cite[Theorem~4.58]{BVoHI}.

\begin{mysubsection}{$\Gamma$-spaces}
\label{sgamma}
Homotopy-coherent abelian monoids may be conveniently described in
terms of a lax version of $\fA$, representing $\Gamma$-spaces 
(cf.\ \cite{SegCC}): 

Let $\Gamma$ denote the category of finite pointed sets, and choose a
set \ $\bn^{+}=\{0,\dotsc,n\}$ \ (with basepoint $0$) for each \ $n\in\NN$. \ 
A \emph{$\Gamma$-object} in a pointed category $\cC$ is a pointed functor \ 
$G:\Gamma\to\cC$; \ the category of all such will be denoted by \
$\Gam{\cC}$. \ Note that if $\cC$ is cocomplete, we can extend $G$ to
all of \ $\Seta$ \ by assuming it commutes with arbitrary colimits.
A $\Gamma$-space $G$ \ -- \ that is, an object in \ $\GS$ \ (or \
$\Gam{\Ta}$) \ -- \ is called \emph{special} if for \
$A,A'\in\Gamma$, \ the natural map \ 
$G(A\vee A')\to G(A)\times G(A')$ \ is a weak equivalence. 
This implies that for each \ $n\in\NN$, \ the obvious map \ 
%
\begin{myeq}\label{eeleven}
G(\bn^{+})\to\underbrace{G(\bo^{+})\times\dotsc\times G(\bo^{+})}_{n}
\end{myeq}
\noindent is a weak equivalence. Such a $G$ is called  \emph{very special} if
in addition \ $\pi_{0}G(\bo^{+})$ \ is an abelian group under the
induced monoid structure. 
\end{mysubsection}

\begin{defn}\label{dclass}
A special $\Gamma$-space $G$ has a \emph{classifying $\Gamma$-space} \ 
$BG$, \ which is itself special, defined by setting \ 
$(BG)(\bn^{+}):=G(\bn^{+}\times\bn^{+})$, \ with the diagonal
structure maps (see \cite[1.3]{SegCC} and compare \cite{MilgB}).
By iterating the functor $B$ we obtain a $\Omega$-spectrum \ 
$\bB G:=\lra{(B^{i}G)(\bo^{+})}_{i=0}^{\infty}$. 
\end{defn}

Thus \ $G(\bo^{+})$ \ itself is an infinite loop space 
(with a specified $H$-space structure) if and only if $G$ is very special.

\begin{mysubsection}{The \ $\Gamma^{+}$-construction}
\label{sqfunct}
For any pointed simplicial set \ $K\in\Sa$, \ Barratt
defines the free simplicial monoid \ $\Gamma^{+}K$ \ to be \
$\coprod_{n\geq 1}~K^{n}\times_{\Sigma_{n}}~W\Sigma_{n}/\sim$, \ where
$\sim$ is generated by the obvious inclusions \ $K^{n}\hra K^{n+1}$ \ 
and \ $\Sigma_{n}\hra\Sigma_{n+1}$ \ (cf.\ \cite[\S 4]{BarF}). \ Then \ 
$\Gamma^{+}K$ \ is actually a $\Gamma$-space (see \cite[\S 8]{AndS}).
To avoid confusion in the notation we shall denote this functor by \ 
$\gamma^{+}:\Sa\to\GS$. \ The (dimensionwise) group completion \ 
$\gamma K:=\Omega B\gamma^{+}K$ \
is a very special $\Gamma $-space, which models the infinite loop
space \ $\Omega^{\infty}\Sigma^{\infty}K$. \ 

The functor \ $\gamma:\Sa\to\GS$ \ is left adjoint to \ 
$G\mapsto G(\bo^{+})$. \ If $K$ is connected, then \ 
$\gamma^{+} K\simeq\gamma K$ \ (cf.\ \cite[Theorem~6.1]{BarF}). 
Note that we can think of \ $\SSa:=\gamma S^{0}$ \ as the inclusion
functor \ $\Gamma\to\Sa$ \ (cf.\ \cite[2.7]{LydaS}).
\end{mysubsection}

\begin{mysubsection}{The model category of $\Gamma$-spaces}
\label{sbfmc}
In \cite[\S 3]{BFrH}, Bousfield and Friedlander define a proper
simplicial model category structure on \ $\GS$ \ as a
diagram category with \ $\Sigma_{n}$-action on each \ $G(\bn^{+})$, \   
which they call the \emph{strict} model category:
a map \ $f:G\to G'$ \ is a weak equivalence if \ 
$f(\bn^{+}):G(\bn^{+})\to G'(\bn^{+})$ \ is a \
$\Sigma_{n}$-equivariant weak equivalence for each \ $n\geq 1$, \  
and it is a (co)fibration if it is a $\Sigma_{n}$-Reedy (co)fibration 
(see \cite[\S 15.3]{PHirM}).

They show that the homotopy category of very special $\Gamma$-spaces is
equivalent to that of connective spectra (see \cite[Theorem~5.1]{BFrH}), 
with Quillen equivalences provided by iterations of the functor $B$
and its adjoint. They then define a \emph{stable} weak equivalence of
$\Gamma$-spaces to be a map inducing a weak equivalence of the
corresponding spectra, and so obtain a new simplicial model category
structure on \ $\GS$ \ (with the same cofibrations, but fewer 
fibrations), whose homotopy category is again equivalent to the usual
stable category of connective spectra (see \cite[Theorem~5.8]{BFrH}). 

Variants on these two model category structures (with the same weak
equivalences) are provided in \cite[App. A]{SchweSH}.
\end{mysubsection}

\begin{mysubsection}{$\Gamma$-simplicial groups}
\label{sgammasg}
In view of \eqref{ethree}, it is natural to think of the category \
$\Gam{\cG}$ \ of $\Gamma$-simplicial groups as representing
connected infinite loop spaces; note that every special
$\Gamma$-object here is trivially \emph{very} special, because of the
shift in indexing for homotopy groups.

A $\Gamma$-simplicial group $G$ also known as a \emph{chain functor}
(cf.\ \cite[\S 1]{AndC}), since one can associate to it a generalized
homology theory by setting \ $H_{n}(X;G):=\pi_{n}(G_{\bullet}X)$ \  
for each \ $X\in\Sa$, \ where the simplicial group \ $G_{\bullet}X$ \ 
is defined by \ $G_{n}X:=G(X_{n})_{n}$. \ 
Here each \ $G(X_{n})\in\cG$ \ is defined as above by extending 
$G$ from $\Gamma$ to \ $\Seta$, \ so that \ $G_{\bullet}X$ \ is
actually the diagonal of a bisimplicial group.

Equivalently, given a $\Gamma$-space \ $G\in\GS$, \
extend it via colimits from $\Gamma$ to\ $\Seta$ \ and thus via the
diagonal to a functor \ $\tilde{G}:\Sa\to\Sa$, \ which in fact 
takes a (pre)spectrum \ $(X_{n})_{n\in\NN}$ \ to a (pre)spectrum \ 
$(\tilde{G}X_{n})_{n\in\NN}$ \ using:
$$
S^{1}\wedge\tilde{G}(X_{n})\to\tilde{G}(S^{1}\wedge X_{n})\to\tilde{G}(X_{n+1})~.
$$
\noindent Thus for each \ $X\in\Sa$, \ one may evaluate the homology
theory associated to $G$ on $X$ by:
$$
H_{n}(X;G)~\cong~\pi_{n}^{S}\tilde{G}(\bS\wedge X)~=~
\colimit{k\to\infty}~\pi_{n+k}\tilde{G}(S^{k}\wedge X)~,
$$
\noindent where \ $\bS:=\lra{S^{n}}_{n=0}^{\infty}$ \ is the sphere spectrum.

Note that if $G$ is very special, then \ $\tilde{G}(\bS\wedge X)$ \ 
is the $\Omega$-spectrum corresponding to Anderson's \
$G_{\bullet}X$ \ (see \cite[\S 4]{BFrH}. 
\end{mysubsection}
 
\begin{mysubsection}{Generalized cohomology}
\label{sgch}
We now explain how the definitions of \S \ref{scohf} apply in this
context: first, note that the usual model category structure on \
$\cE=\Sa$ \ is symmetric monoidal and enriched over \ $\cV=\Sa$ \
(cf.\ \cite[II, \S 3]{QuiH}). Now for \ $\Phi=\Gamma$, \ Lydakis 
(in \cite{LydaS}) defined a smash product of $\Gamma$-spaces making \
$\PhV=\GS$, \ too, into a symmetric monoidal category, with unit
$\SSa$. He also defines internal function complexes \
$\uHom_{\GS}(G,H)\in\GS$ \ for \ $G,H\in\GS$ \ by setting: 
%
\begin{myeq}\label{enine}
\uHom_{\GS}(G,H)(\bn^{+})~:=~\map_{\GS}(G,H(\bn^{+}\wedge-))~,
\end{myeq}
\noindent 
where \ $H(\bn^{+}\wedge-))(\bk^{+}):=H(\bn^{+}\wedge\bk^{+})$ \ and \
$\map_{\GS}(-,-)\in\Sa$ \ is the usual simplicial function complex. 

Thus \ $\PhE=\GS$ \ is indeed enriched over \ $\PhV$ \ (cf.\ 
\cite[2.1]{LydaS}). Moreover, \ $\PhV$ \ is semi-triangulated, with
the delooping \ $B:\GS\to\GS$ \ (\S \ref{dclass}) as the ``suspension
automorphism'' $T$ of \S \ref{esgps}. The deloopings of the $0$-sphere \ 
$\{B^{n}\SSa\}_{n=0}^{\infty}$ \ corepresent homotopy groups in \
$\ho\GS$, \ since its homotopy category is equivalent to that of
connective spectra, with generator $\SSa$ (corresponding to \ $S^{0}$). 

Now for any $\Gamma$-space \ $G\in\PhE$ \ and any pointed simplicial 
set \ $K\in\cE$, \ $\Hom\h_{\cE}(K,G)$ \ is a fibrant $\Gamma$-set 
(\S \ref{rhom}), so the $\Sa$-function complex \ $M:=\maps(K,G)$ \ is
a $\Gamma$-space. If $G$ is (very) special, so is $M$, since \
$\maps(K,-)$ \ has homotopy meaning and preserves products.

Moreover, applying Barratt's functor yields a special $\Gamma$-space \
$\gamma K$, \ and the adjunction isomorphism: 
%
\begin{myeq}\label{eten}
M=\maps(K,G)~\xra{\cong}~\uHom_{\PhE}(\gamma K,G)
\end{myeq}
\noindent induces an isomorphism between the homotopy groups of $M$
and those of \ $\uHom{\PhE}(\gamma K,G)$ \ (corepresented by $\SSa$ and
its suspensions). 

Therefore, for special $G$ the homotopy groups of $M$ are determined by
those of \ $M(\bo^{+})=\maps(K,G(\bo^{+}))$, \ which are by definition \ 
$H^{\ast}(K;G)$, \ the generalized cohomology groups associated
to the $\Omega$-spectrum for $G$.
\end{mysubsection}
 
\begin{mysubsection}{Generalized homology}
\label{sgh}
Barratt's functor \ $\gamma:\cE\to\PhE$ \ is the required
functor \ $\Ap$, \ by \ \eqref{eten}, \ so its left derived functors \ 
are \ $\pi_{\ast}\gamma K$ \ (since every $K$ is cofibrant). These 
turn out to be the stable homotopy groups of $K$, and are by
definition the homology groups of $K$ in this context.

Finally, since the smash product of (cofibrant) $\Gamma$-spaces is
taken to the smash product of spectra under the equivalence of
homotopy categories (see \cite[Lemma 5.16]{LydaS}), we see that the
groups \ $H_{\ast}(K;G)$ \ of \S \ref{shomol} are just the generalized 
homology groups associated to the $\Omega$-spectrum for $G$.
\end{mysubsection}

\begin{mysubsect}{The (co)simplicial version}
\label{scsv}

We next show how these definitions can be made to fit the description
in \S \ref{scver}\vsn : 

First, note that \ $s\cS$, \ as well as \ $s\Ta$ \ and \ $s\cG$ \ 
(cf.\ \S \ref{rtops}), \ have resolution model category structures
with \ $\cM=\{S^{1}\}$ \ -- \ this is the original \ $E^{2}$-model
category of \cite[\S 5.10]{DKStoE}, \ which was constructed precisely
so that if \ $\Vd$ \ is a resolution of \ $X\in\cS$, \
then the diagional \ $\diag \Vd$ \ (or equivalently, the realization of
the corresponding simplicial space) is weakly equivalent to $X$.
Moreover, $\cS$, as well as \ $\Ta$ \ and $\cG$, are enriched over \
$\cV:=\cS$ \ with its usual closed symmetric monoidal structure.

We also need a suitable model category structure on the category \
$c\GS$ \ of cosimplicial $\Gamma$-spaces \ -- \  namely, the dual
of Moerdijk's model category of bisimplicial sets (cf.\ \cite[\S 1]{MoerB}), 
in which a map \ $f:\Xu\to\Yu$ \ of cosimplicial $\Gamma$-spaces is a 
weak equivalence (resp., cofibration) \ if \ $\Tot f$ \ is a
weak equivalence (resp., cofibration) of $\Gamma$-spaces. This implies
that \ $\Tot:c\GS\to\GS$ \ induces an equivalence of
homotopy categories, so for all practical purposes we can avoid
working with cosimplicial objects altogether (but see Theorem
\ref{tthree} below). The inverse equivalence \ $c\GS\to c\GS$ \ is defined by
$\Phi\mapsto\cu{\Phi}$ \ (the constant cosimplicial object).
Thus \ $\ho(c\GS)$ \ (with this structure) is equivalent to the 
stable category of connective spectra, which is semi-triangulated,
with \ $\cu{B}\circ\Tot:c\GS\to c\GS$ \ (\S \ref{dclass}) 
as the suspension automorphism $T$, and \ $\cu{\SSa}$ \ as generator.

Now, given a special $\Gamma$-space \ $G\in\GS$ \ and a
free simplicial resolution \ $\Vd\to X$ \ in the original resolution
model category \ $s\cS$, \ for any simplicial set $Y$ \ -- \ 
in particular, for \ $Y=G(\bo^{+})$ \ -- \ we have:
%
\begin{myeq}\label{etwentysix}
\maps(\diag\Vd,Y)\cong\Tot\maps(\Vd,Y)
\end{myeq}
\noindent (see \cite[XII, \S 4.3]{BKaH}). Thus in our case the
cosimplicial $\Gamma$-space \ $\maps(\Vd,G)$ \ 
is weakly equivalent to the (constant cosimplicial) space \ 
$\cu{\maps(X,G(\bo^{+}))}$, \ whose homotopy groups are \ $H^{\ast}(K,G)$ \ 
(\S \ref{sgch}).

Finally, note that Barratt's functors \ $\gamma^{+}$ \ and $\gamma$
are defined dimensionwise on a simplicial set $K$, \ so that \ 
$$
\diag\gamma\Vd~=~\gamma\diag\Vd
$$
for any bisimplicial set \ $\Vd$. \ 
Thus we may define \ $\Ap:\cE\to\PhE$ \ to be $\gamma$, \ and its total left
derived functor is naturally equivalent to $\gamma$ (in Moerdijk's
model category \ $s\Sa$), \ since \ $\diag\Vd\xra{\simeq}K$ \ for any
free simplicial resolution \ $\Vd\to K$. \ Thus again the (unadorned)
homology groups are the stable homotopy groups of $K$, and \
$H_{\ast}(K;G)$ \ are the generalized homology groups associated to
the $\Omega$-spectrum for $G$. 
\end{mysubsect}

%
%
\section{The spectral sequences}
\label{css}

We now want to use this machinery  to try to understand relationships
among the various homology and cohomology theories. First, we shall
need a preliminary notion: 

\begin{defn}\label{dcpalg}
If $\cM$ is a set of models in a model category $\cC$ \ (with \ 
$\PiM\subseteq\cC$ \ as in \S \ref{srmc}), then \ $\CP:=(\ho\PiM)\op$ \ is 
a $\fG$-theory, which sketches the category \ $\CPA$ \ of
\emph{$\CP$-algebras} (cf.\ \cite[\S 3]{BStoG}).
\end{defn}

\begin{remark}\label{rcpalg}
If we think of $\cM$ and its suspensions as corepresenting homotopy
groups in $\cC$ \ (cf.\ \S \ref{rsphx}), \ then \ $\CP$-algebras are 
graded groups equipped with an action of the corresponding primary
homotopy operations \ - \ the motivating example being \
$\picM{\ast}X$ \ for any \ $X\in\sC$. \
This notion may be extended to any concrete category $\cC$ by the
conventions of \cite[\S 3.2.2]{BStoG}, and may also be dualized as in
\cite{BouC} by taking \ $\CP:=\ho\PiM$, \ rather than the opposite
category (cf.\ \cite[\S 1.13]{BPescF}).  

Note that the derived functors of any functor into $\cC$ actually
take values in \ $\CPA$. 
\end{remark}

\begin{examples}\label{ecpalg}
\begin{enumerate}\renewcommand{\labelenumi}{(\alph{enumi})}
\item If $\cC$ has a trivial model category structure, and $\cM$
  consists of (enough) projective generators \ -- \ e.g., if
  $\cC=\TA$ \ and \ $\cM=\fFTp$  \ -- \ then \ $\CPA\cong\cC$.
\item If \ $\cC=s\cDD$ \ or \ $c\cDD$ \ for some abelian category \
  $\cDD$, \ and $\cM$ again consists of (enough) projective 
  generators \ -- \ e.g., for \ $\cC=s\TX$ \ and $\cM$ as above \ -- \
  then \ $\CPA\cong\gr_{\NN}\cDD$ \ (where we use lower or upper
  indices for the grading according to the usual convention).
\item For \ $\cC=\Ta$ \ or \ $\Sa$, \ with \ $\cM=\{S^{1}\}$, \
  then \ $\CPA\cong\PAlg$ \ is the category of ordinary
  $\Pi$-algebras, modeling the usual homotopy groups of topological
  spaces.
\item If \ $\cC=\GS$ \ and \ $\cM=\{\SSa\}$, \ then \
  $\CPA$ \ is equivalent to the category of graded connected 
  $\bp$-modules for \ $\bp=\pi^{S}_{\ast}S^{0}$ \ (homotopy groups of
  the sphere spectrum), since \ $\piM{\ast}G$ \ are just the 
  stable homotopy groups of the $\Omega$-spectrum corresponding to \
  $G\in\GS$.
\end{enumerate}

Using the Quillen equivalence of \ \eqref{ethree}, \ we see that when \ 
$\cC=s\TA$ \ we often have interesting categories of $\CP$-algebras
(see, e.g., \cite[\S 3.2.1]{BStoG}). 
\end{examples}

We shall also need the following version of \cite[Prop.~3.2.3]{BStoG}:

\begin{prop}\label{pind}
Any contravariant functor \  $T:\cC\to c\cB$ \ from a model
category $\cC$ (equipped with a set of models $\cM$) to a concrete
category $\cB$  induces a graded functor \  
$\bar{T}^{\ast}:s\CPA\to s\BPA$ \ by setting \  
$\bar{T}^{k}(\picM{\ast}\Vd):=\pi^{k}(T\Vd)$ \ for cofibrant \
$\Vd\in\sC$, \ and extending by taking $0$-th derived functor.
\end{prop}

\begin{proof}
Since \ $\picM{\ast}:\ho\PiM\to\fF_{\CP}$ \ is an equivalence of 
categories (onto the free \ $\CP$-algebras), in particular \ 
$\picM{\ast}\Vd\cong\picM{\ast}\Wd~\EQUIV~\Vd\simeq\Wd$ \ for cofibrant \ 
$\Vd,\Wd\in\sC$, \ so \ $\bar{T}^{\ast}$ \ is well-defined on free \
$\CP$-algebras. 
\end{proof}

\begin{mysubsect}{A general setting}
\label{sgen}

In Sections \ref{caft}-\ref{cgencoh} the algebraic and topological
versions of homology and cohomology have been treated separately. We now
show how the Procrustean framework of \S \ref{scver} may be used in
order to obtain a uniform description of various relations between them. 
\end{mysubsect}

\begin{examples}\label{egset}
We wish to concentrate on the following list of cohomological settings 
(Definition \ref{dcset}), discussed above:
\begin{enumerate}\renewcommand{\labelenumi}{(\alph{enumi})~}
\item $\lra{\cC=\TA,\cM=\fFTp,\cV=\Set,\Phi=\fA,\Ap=\Aa}$ \ 
      for some $\fG$-theory $\Theta$;
\item More generally, \ 
      $\lra{\cC=\TA/X,\cM=\fFTp/X,\cV=\Set,\Phi=\fA,\Ap=\Ax}$ \ 
      for some $\fG$-theory $\Theta$ and fixed \ $X\in\TA$.
\item $\lra{\cC=s\TA/X,\cM=\{\co{\FT(s)}~|\ \FT(s)\in\fFTp\},
      \cV=\Sa,\Phi,\Ap}$ \ where $\Phi$ is some strong $\fA$-sketch.
\item $\lra{\cC=\Sa,\cM=\{S^{1}\},\cV=\Sa,\Phi=\Gamma,\Ap=\gamma}$ \ 
      (with the symmetric monoidal structure on \ $\GS$ \ of
      \S \ref{sgch}). 
\end{enumerate}
\end{examples}

In all these examples we have additional properties which
we shall require in our applications, which we may formalize as follows:

\begin{defn}\label{dcomps}
A cohomological setting \ $\lra{\cC,\cM,\cV,\Phi,\Ap}$ \ is 
\emph{complete} if if it is equipped with: 

\begin{enumerate}
\item A left adjoint \ $\diag:\sC\to\cC$ \ to the inclusion \ 
  $\co{-}:\cC\to\sC$, \ which induces \ $\diag:s\PhC\to\PhC$, \ as well
  as a convergent first-quadrant spectral sequence with:
%
\begin{myeq}\label{etwentythree}
E^{2}_{s,t}~\cong~\pi_{s}\piM{t}\Vd~\Longrightarrow~\piM{s+t}(\diag\Vd)~,
\end{myeq}
\noindent for each \ $\Vd\in\sC$ \ and \ $M\in\cM$;
\item A right adjoint \ $\Tot:c\cV\to\cV$ \ to the inclusion \ 
  $\cu{-}:\cV\to c\cV$, \ which induces \ $\Tot:c\PhV\to\PhV$, \ as
  well as a second-quadrant spectral sequence with:
%
\begin{myeq}\label{etwentyfour}
E_{2}^{s,t}~\cong~\pi^{s}\pin{M'}{t}\Xu~\Longrightarrow~
\pin{M'}{t-s}(\Tot\Xu)~,
\end{myeq}
\noindent for each \ $\Xu\in c\cV$ \ and \ $M'\in\MP$ \ (we do not 
address questions of convergence);
\item A natural ``$\PhC$-adjointness'' isomorphism: 
%
\begin{myeq}\label{etwentyfive}
\Tot(\uHom(\Vd,G)~\xra{\cong}~\uHom(\diag\Vd,G)
\end{myeq}
\noindent for any \ $\Vd\in s\PhC$ \ and \ $G\in\PhC$.
\end{enumerate}
\end{defn}

\begin{prop}\label{pccs}
Each of the examples of \S \ref{egset} is a complete cohomological
setting.
\end{prop}

\begin{proof}
Since (a) and (b) are instances of (c), we have only two cases to
consider\vsm:

\noindent(1)\  Assume \ $\cC=s\TA/X$ \ for some $\fG$-theory $\Theta$
  sorted by $S$. Then \ $\Vd\in\sC$ \ is a bisimplicial \Tal\ (over $X$), 
  and let \ $\diag\Vd$ \ be the usual diagonal (with \ 
  $(\diag\Vd)_{n}:=(V_{n})_{n}$). \ 
  Note that \ $\UT\Vd$ \ is just an an $S$-graded bisimplicial set,
  with \ $\UT\diag\Vd=\diag\UT\Vd$ \ (even though colimits are not
  generally preserved by \ $\UT$). \  By Remark \ref{rsphx} we see
  that the Bousfield-Friedlander spectral sequence for \ $\UT\Vd$ \ in
  each degree (cf.\ \cite[Theorem~B.5]{BFrH}) has the form \
  \eqref{etwentythree}. 

  Similarly, given a cosimplicial object \ $\Xu\in c(s\PhTAX)$, \ 
  the usual \ $\Tot$ \ for the ($S$-graded) cosimplicial simplicial set \ 
  $\UT\Xu$ \ is defined to be the simplicial set \ $T_{\bullet}$ \ with \ 
  $T_{n}:=\Hom_{c\Set}(\bDb\otimes\Delta[n],\Xu)$, \ 
  and this has a natural structure of a \Phal\ in \ $\TA/X$ \ by
  Remarks \ref{rhom} and \ref{rcomma} and \S \ref{sspat}.
  Thus \ $\Tot\UT\Xu$ \ lifts to \ $\Tot\Xu\in s\TA$. \ The homotopy 
  spectral sequence for the cosimplicial space \ $\UT\Xu$, \ with:
$$
E_{2}^{s,t}~=~\pi^{s}\pi_{t}\UT\Xu~\Longrightarrow~\pi_{t-s}(\Tot\UT\Xu)~,
$$
  (see \cite[X, 6.1 \& 7.2]{BKaH}) \ gives \ \eqref{etwentyfour} \ 
  (though it does not necessarily converge!).

  Finally, \ \eqref{etwentyfive} \ follows from \ \eqref{etwentysix}\vsm.

\noindent (2)\  For \ $\cC=\Sa$ \ we can use the usual diagonal and \ $\Tot$ \ 
  and the original spectral sequences for (co)simplicial spaces. 
  For \ \eqref{etwentyfive}, \ consider the cosimplicial
  $\Gamma$-space \ $\Eu:=\uHom_{\GS}(\Vd,G)$: \ Definition \
  \eqref{enine} \ of \ $\uHom_{\GS}$ \ in terms of the simplicial
  function complex \ $\map_{\GS}$ \ shows that \ 
  $\Tot\Eu~\cong~\uHom_{\GS}(\diag\Vd,G)$ \ again, by \ \eqref{etwentysix}.
\end{proof}

With this at hand, we can describe several spectral sequences
connecting the various functors we have defined so far. First, a
universal coefficients theorem for cohomology:

\begin{thm}\label{tone}
Let \  \ $\lra{\cC,\cM,\cV,\Phi,\Ap}$ \ be a complete cohomological
setting, and let $G$ be a \Phal\ in $\cC$. Then for any \ $Y\in\cC$ \
there is a natural cohomological spectral sequence with
$$
E_{2}^{s,t}~\cong~\bExt^{s,t}(H_{\ast}Y,G)~
\Longrightarrow~H^{t-s}(Y;G)~,
$$
\noindent where \ $\bExt^{s,t}(C,G):=(L_{s}\bar{T}(C))^{t}$ \ for any \
$C\in\hy{(\PhC)}{\PAlg}$, \ and \ $T:=\uHom(-,G)$. 
\end{thm}

\begin{proof}
Let \  $Z\to Y$ be a cofibrant replacement in \ $\cC$, \ and assume
$G$ is fibrant. 
We use \ $\MP:=\{\Ap M\}_{M\in\cM}$ \ as models in \ $\PhC$ \ (\S
\ref{scver}), with \ $T^{n}$ \ as the suspension (\S \ref{striang}),
to define the resolution model category structure on \ $s\PhC$. \ 
As in the proof of \cite[Theorem~4.2]{BStoG}, let \ $\Vd\to\Ap Z$ \ be
a free simplicial resolution in  \ $s\PhC$, \ so that by \
\eqref{etwentythree} \ the natural map \ $\diag\Vd\to\Ap Z$ \ is a
weak equivalence. 

If we set \ $\Eu:=\uHom(\Vd,G)$ \ (a cosimplicial \Phal\ in $\cC$),
then by \ \eqref{etwentyfive} \ and \  \ \eqref{eone}:
$$
\Tot\Eu~=~\uHom(\diag\Vd,G)~\simeq~\uHom(\Ap Z,G)~\cong\map(Z,G)
~=~\LL\map(-,G)(Y)
$$
ao \ $\pin{\MP}{t-s}(\Tot\Eu)=\pin{\MP}{t-s}\map(Z,G)=H^{t-s}(Y;G)$ \ 
by Definition \ref{dscoh}.

On the other hand, since each \ $V_{n}$ \ is cofibrant:
$$
\pin{\MP}{\ast}E^{n}=\pin{\MP}{\ast}\uHom(V_{n},G)~=~
\bar{T}(\pin{\MP}{\ast}V_{n})
$$
and since \ $\Vd\to\Ap Z$ \ is a cofibrant replacement, \ 
$\pin{\MP}{\ast}\Vd\to\pin{\MP}{\ast}\Ap Z=:H_{\ast}Y$ \ is a 
free resolution in \ $\hy{(\PhC)}{\PAlg}$, \ so:
$$
\pi^{s}\pin{\MP}{\ast}\Eu~=~
\pi^{s}(\bar{T}(\pin{\MP}{\ast}\Vd))~=~
\pi^{s}\LL\bar{T}(H_{\ast}Y)~=~
L^{s}\bar{T}(H_{\ast}Y)~,
$$
as claimed.
\end{proof}

Note that for generalized cohomology of spaces this takes the 
familar form (cf.\ \cite{AdG} and \cite[IV, \S 4]{EKMMayR}):

\begin{cor}\label{cone}
For any special \ $G\in\GS$ \ and \ $K\in\Sa$ \ 
there is a second quadrant spectral sequence with:
$$
E_{2}^{s,\ast}~\cong~\Ext_{\RM{\bp}}^{s}(\pi^{S}_{\ast}K,G)~
\Longrightarrow~H^{s-t}(K;G)\vsm.
$$
\end{cor}

There is also a version for homology:

\begin{prop}\label{ptwo}
Let \  \ $\lra{\cC,\cM,\cV,\Phi,\Ap}$ \ be a complete cohomological
setting, and let $G$ be a \Phal\ in $\cC$. Then for any \ $Y\in\cC$ \
there is a natural first quadrant spectral sequence with
%
\begin{myeq}\label{etwentyeight}
E^{2}_{s,t}~\cong~\bTor_{s,t}(H_{\ast}Y,G)~
\Longrightarrow~H_{t+s}(Y;G)~,
\end{myeq}
\noindent where \ $\bTor_{s,\ast}(C,G):=(L_{s}\bar{T}(C))$ \ for any \
$C\in\hy{(\PhC)}{\PAlg}$, \ and \ $T:=-\otimes G$. 
\end{prop}

\begin{proof}
This generalization of \cite[Theorem~4.4]{BStoG} for the composite functor:
$$
\PiM~\xra{\Ap}~\PhC~\xra{-\otimes G}~\PhC
$$
is proven like Theorem \ref{tone}, with \ \eqref{etwentythree} \ 
replacing \ \eqref{etwentyfour}\vsm .
\end{proof}

For generalized homology this takes the form:

\begin{cor}\label{ctwo}
For any special \ $G\in\GS$ \ and \ $K\in\Sa$ \ 
there is a natural first quadrant spectral sequence with:
$$
E^{2}_{s,t}~\cong~\Tor^{\RM{\bp}}_{s,t}(\pi^{S}_{\ast}K,G)~
\Longrightarrow~H_{t+s}(K;G)\vsm.
$$
\end{cor}

Finally, we have the following two generalizations of \cite{BlaH}:

\begin{thm}\label{ttwo}
Let \  \ $\lra{\cC,\cM,\cV,\Phi,\Ap}$ \ be a complete cohomological
setting, and let $G$ be a \Phal\ in $\cC$. Then for any \ $Y\in\cC$ \
there is a natural first quadrant spectral sequence with
$$
E^{2}_{s,t}~\cong~ L_{s}\bar{T}(\picM{\ast}Y)_{t}~\Longrightarrow~H_{t+s}(Y;G)~,
$$
\noindent where \ where \ $T:=\Ap(-)\otimes G$.
\end{thm}

\begin{proof}
Similar to the proof of Theorem \ref{tone}, except that here we start
with a free simplicial resolution \ $\Vd\to Y$ \ in \ $\sC$, \ and
note that in this case \ $\picM{\ast}\Vd\to\picM{\ast}Y$ \ 
is a free simplicial resolution in the category \ $s\hy{\cC}{\PAlg}$.
\end{proof}

In \cite[Prop.~5.1]{SegCS}, Segal produced a stable version of this 
spectral sequence for any generalized homology theory \ $k_{\ast}$ \
(converging strongly to \ $k_{\ast}X$ \ if\ \ $k_{\ast}$\ \ is
connective).

\begin{thm}\label{tthree}
For $Y$ and $G$ as above, there is a natural second quadrant spectral
sequence with: 
$$
E_{2}^{s,t}~\cong~\wExt^{s}_{t}(\picM{\ast}Y,G)~\Longrightarrow~H^{t-s}(X;G)~,
$$
\noindent where \ $\wExt^{s}(-,G):=L_{s}\bar{T}$ \ for \ 
$T:=\map_{\cC}(-,G)$.
\end{thm}

Note that Schwede, in \cite[\S 5.5]{SchweSA}, also defined a spectral
sequence relating the stable homotopy of a \Tal\ to Quillen homology.

%
%

%
\end{document}